%% file: cwnew.tex
\documentclass[a4paper,10pt,nointlimits]{article}

\input{preamble.tex}

\begin{document}
\title{Extending Chacon-Walsh: Minimality and Generalised Starting Distributions}
\author{A.~M.~G.~Cox\thanks{The author wishes to thank David Hobson for his helpful suggestions and advice. e-mail:
        \texttt{amgc500@york.ac.uk}; web:
        \texttt{www-users.york.ac.uk/$\sim$amgc500/}}\\
        Department of Mathematics,\\
        University of York,\\
        York Y010 5DD, U.~K.}
% \and D.~G.~Hobson\thanks{e-mail:
%        \texttt{dgh@maths.bath.ac.uk}; web:
%        \texttt{www.bath.ac.uk/$\sim$masdgh/}},\\
%        Department of Mathematical Sciences,\\
%        University of Bath,\\
%        Bath BA2 7AY, U.~K.}
\maketitle

\begin{abstract}
In this paper
we consider the Skorokhod embedding problem for general starting and target measures. In particular, we provide necessary and sufficient conditions for a stopping time to be minimal in the sense of \citet{Monroe:72}. The resulting conditions have a nice interpretation in the graphical picture of \citet{ChaconWalsh:76}.

Further, we demonstrate how the construction of Chacon and Walsh can be extended to any (integrable) starting and target distributions, allowing the constructions of %\citet{AzemaYor:79,Jacka:88,Vallois:83}
Azema-Yor, Vallois and Jacka
to be viewed in this context, and thus extended easily to general starting and target distributions. In particular, we describe in detail the extension of the Azema-Yor embedding in this context, and show that it retains its optimality property.
\end{abstract}

\input{introcw.tex}

\input{basicbalayage.tex}

%\input{newbalayage.tex}

\input{minbasic.tex}

\input{minandpot.tex}

\input{minlimit.tex}

\input{AYtype.tex}

\bibliographystyle{jmr}
\bibliography{general}

\end{document}

%% file: preamble.tex
\usepackage{amssymb,amsmath,amsthm,verbatim}
\usepackage[dvips,final]{graphicx}
\usepackage{psfrag}
\usepackage{natbib}
\newcommand{\R}{\mathbb{R}}
\newcommand{\Prob}{\mathbb{P}}
\renewcommand{\Pr}{\Prob}
\newcommand{\supp}{\mbox{supp}}

\newcommand{\E}{\mathbb{E}}

\newcommand{\Fc}{\mathcal{F}}
\newcommand{\Lc}{\mathcal{L}}

\newcommand{\scA}{\mathcal{A}}
\newcommand{\Borel}{\mathcal{B}}
\newcommand{\half}{\frac{1}{2}}
\newcommand{\N}{\mathbb{N}}
\newcommand{\Z}{\mathbb{Z}}
\newcommand{\eps}{\varepsilon}

\newtheorem{theorem}{Theorem}[section]
\newtheorem{prop}[theorem]{Proposition}
\newtheorem{lemma}[theorem]{Lemma}

\theoremstyle{definition}

\newtheorem{remark}[theorem]{Remark}
\newtheorem{definition}[theorem]{Definition}

\newcommand{\indic}[1]{\boldsymbol{1}_{\{\ensuremath{#1}\}}}

\newcommand{\ul}[1]{\underline{#1}}
\newcommand{\ol}[1]{\overline{#1}}

\DeclareMathOperator{\sign}{sign}

\DeclareMathOperator*{\argmin}{argmin}
\newcommand{\as}{a.s.}
\newcommand{\ie}{i.e.}

%% file: introcw.tex
%In this chapter we consider a more general method for constructing stopping times, following the method of \citet{ChaconWalsh:76}. In this context we are able to consider the embedding problem where the process has an integrable starting distribution. Consideration of general starting distributions requires a more general characterisation of minimality, and we are able to use this characterisation of minimality to provide a simple condition for the construction to be minimal.
%
%We are then able to use the construction to provide a simple description of the stopping times considered in Chapter~\ref{ch:minAY} --- in particular, the Azema-Yor and modulus maximising stopping times can be easily extended to the case where we have a more general starting distribution, and shown to be minimal and optimal. Further, we are able to show that the stopping time introduced by \citet{Vallois:83} can be described simply in the construction, which allows us to extend the stopping time to general starting measures.
%
%
%\newpage

\section{Introduction}
The Skorokhod embedding problem has a long history, and was first posed (and solved) in \citet{Skorokhod:65}. Simply stated it is the following: given a stochastic process $(X_t)_{t \ge 0}$ and a distribution $\mu$, find a stopping time $T$ such that $X_T \sim \mu$.

In this work we will be interested in the case where $B_t$ is a Brownian motion on $\R$, with a given (integrable) starting distribution $\mu_0$. Since Brownian motion on $\R$ is recurrent, the existence of such a stopping time is trivial: consider an independent random variable $Y$ with distribution $\mu$ and run until the first time that the Brownian motion hits $Y$. Hence interest lies in the properties of the stopping time $T$ and also of the stopped process $B_{t \wedge T}$.

Classically, the 1-dimensional question has been considered in the case where $B_0 = 0$, and the target distribution $\mu$ is centred. In this case many solutions are known: \citet{AzemaYor:79,Bertoin:92,Dubins:68,Jacka:88,Perkins:86,Root:69,Vallois:83}. We refer the reader to \citet{Obloj:04b} for an excellent recent survey of these results. A property shared by all of these embeddings is that the process $B_{t \wedge T}$ is uniformly integrable, and we shall call stopping times for which this is the case UI stopping times. Further, within the class of embeddings where $T$ is UI, many of these stopping times have optimality properties: for example the Azema-Yor embedding maximises the law of the maximum, while the Vallois construction can be used to minimise or maximise $\E(f(L_T))$ for a convex function $f$ \citep{Vallois:92}. It is clear that either of the maximisation problems are degenerate when looked at outside this class.

The class of UI stopping times can also be characterised in the following way due to \citet{Monroe:72}. We make the following definition:
\begin{definition} \label{def:minimal}
A stopping time $T$ for the process $X$ is {\em minimal} if
whenever $S \le  T$ is a stopping time such that $X_S$ and $X_T$
have the same distribution then $S = T$ \as.
\end{definition}
Then the class of minimal stopping times can be shown to be equivalent to the class of UI embeddings we had before:
\begin{theorem}{\citep[Theorem 3]{Monroe:72}}\label{thm:monroeB}
Let $S$ be a stopping time such that $\E(B_S) = 0$. Then $S$ is
minimal if and only if the process $B_{t \wedge S}$ is uniformly
integrable.
\end{theorem}
Such a characterisation gives a natural interpretation to the class of UI embeddings.

Our interest in this paper lies in the extension to general starting measures. In such an example, even if the means agree, there is no guarantee that there will exist a UI stopping time which has the given starting and target distributions. This can be seen by considering the example of a target distribution consisting of a point mass at zero, but with starting distribution of mass $\half$ at each of $-1$ and $1$. Clearly the only minimal stopping time is to stop the first time the process hits $0$, however this stopping time is not UI.

In \citet{CoxHobson:03} conditions for a stopping time to be minimal were considered. When the Brownian motion starts at the origin, and the target distribution is not centred, conditions on the process can be given which are equivalent to the stopping time being minimal. One of the main results of this work is to show that the conditions have an extension to the case of a general starting distribution, however the simple example given above shows that the extension is not trivial.

It will turn out that the characterisation of minimal stopping times is closely connected to the potentials of the two measures. In this context, the relationship between the measures can be viewed graphically in the framework of \citet{ChaconWalsh:76}. In this paper a graphical construction is interpreted as a sequence of exit times from compact intervals, whose limit is an embedding. This is done for starting and target distributions which satisfy the relationship
\begin{equation} \label{eqn:introcond}
- \E^{\mu_0} |X-x| = u_{\mu_0}(x) \ge u_{\mu}(x) = - \E^{\mu} |X-x|
\end{equation}
for all $x \in \R$. We shall show that the construction can be extended to the case where this condition fails, and that the exact method of the extension will determine whether the stopping time is minimal.

Establishing this connection will then allow us to extend several existing embeddings \citep{AzemaYor:79,Jacka:88,Vallois:92}, to the more general setting (maintaining minimality) via a limiting argument.

%% file: basicbalayage.tex
\section{The Balayage Construction} \label{sec:basicbalayage}

In the theory of general Markov processes, a common definition of the
potential of a stochastic process is given by
\begin{equation*}
U \mu (x) = \int_\R \mu(dy)\, \int_{\R^+} ds \, p_s(x,y),
\end{equation*}
where $p_s(x,\cdot)$ is the transition density at time $s$ of the process started
at $x$. In the case of Brownian motion, we note that the
integral is infinite. To resolve this we use the compensated
definition (and introduce new notation to emphasise the fact that this
is not the classical definition of potential):
\begin{equation*}
u_{ \mu} (x) = \int_\R \mu(dy) \,\int_{\R^+} ds \,(p_s(x,y)-p_s(0,0)).
\end{equation*}
This definition simplifies to the following:
\begin{equation} \label{eqn:potdefn}
u_\mu(x) = -\int |x-y| \, \mu(dy).
\end{equation}

\begin{remark}\label{rem:potprop}
The function $u_{\mu}$ has the following properties:
\begin{enumerate}
\item
The measure $\mu$ is integrable if and only if the function $u_\mu$ is finite
for any (and therefore all) $x \in \R$.
\item
$u_{\mu}$ is continuous, differentiable everywhere
except the set $\{x \in \R: \mu(\{x\}) > 0\}$ and concave.
\item
Write
\begin{equation*}
m = \int x \, \mu(dx).
\end{equation*}
As $|x| \to \infty$, we have
\begin{equation} \label{eqn:pmulimit}
u_\mu(x) +|x|  \to
m \sign(x).
\end{equation}

\item
As a consequence of \eqref{eqn:pmulimit}, if $\mu$ and
$\nu$ are integrable distributions, then there exists a constant $K >
0$ such that:
\begin{equation*}
\sup_{x \in \R} |u_{\mu}(x) - u_\nu(x)| <K.
\end{equation*}
\item
$u_{\mu}$ is almost everywhere differentiable with left
and right derivatives
\begin{eqnarray*}
u_{\mu,-}'(x) & = & 1 - 2 \mu((-\infty,x));\\
u_{\mu,+}'(x) & = & 1 - 2 \mu((-\infty,x]).
\end{eqnarray*}
\end{enumerate}
\end{remark}

\citet{Chacon:77} contains many results concerning these potentials. We will
describe a balayage technique that produces a sequence of measures and
corresponding stopping times, and which will have as its limit our
desired embedding. The following lemma will therefore be important in
concluding that the limit we obtain will indeed be the desired
distribution:
\begin{lemma}[\citet{Chacon:77}, Lemmas~2.5, 2.6]\label{lem:limexists}
Suppose $\{\mu_n\}$ is a sequence of probability measures. If
\begin{enumerate}
\item
$\mu_n$ converges weakly to $\mu$ and $\lim_{n \to \infty}
u_{\mu_n}(x_0)$ exists for some $x_0 \in \R$, then \linebreak $\lim_{n \to
\infty} u_{\mu_n}(x)$ exists for all $x \in \R$ and there exists $C
\ge 0$ such that
\begin{equation} \label{eqn:munlimit}
\lim_{n \to \infty} u_{\mu_n}(x) = u_{\mu}(x) -C.
\end{equation}

\item
$\lim_{n \to \infty} u_{\mu_n}(x)$ exists for all $x \in \R$ then
$\mu_n$ converges weakly to $\mu$ for some measure $\mu$ and $\mu$ is
uniquely determined by the limit $\lim_{n} u_{\mu_n}(x)$.
\end{enumerate} 
\end{lemma}
%\psfrag{lambda1}{$\lambda_1$}
%\psfrag{-gp}{$-\gamma_+(\lambda_1)$}
%\psfrag{lambda2}{-$\lambda_2$}
%\psfrag{gm}{$\gamma_-(\lambda_2)$}
%\psfrag{c(x)}{$\kappa(x)$}
%\psfrag{x}{$x$}
%\psfrag{int}{$\beta$}
%\psfrag{alpha}{$\alpha$}
%\psfrag{St}{$\ol{M}_t$}
%\psfrag{It}{$-\ul{M}_t$}
%\psfrag{g(I)}{$\gamma_-(-\ul{M}_t)$}
%\psfrag{g(S)}{$\gamma_+(\ol{M}_t)$}
%\psfrag{gn(I)}{$\gamma_-^{\mu^n}(-\ul{M}_t)$}
%\psfrag{gn(S)}{$\gamma_+^{\mu^n}(\ol{M}_t)$}
%%\psfrag{n1}{$n \vee (-\xi_-)$}
%%\psfrag{n2}{$n \vee \xi_+$}
%\psfrag{n1}{$-\xi_-$}
%\psfrag{n2}{$\xi_+$}
%\psfrag{n3}{$n$}
%\psfrag{n4}{$n$}
%\psfrag{gm(I)}{$\gamma_-^{\mu}(-\ul{M}_t)$}
%\psfrag{gm(S)}{$\gamma_+^{\mu}(\ol{M}_t)$}
%\psfrag{e}{$\eps$}
%\psfrag{a1}{$a_1$}
%\psfrag{a2}{$a_2$}
%\psfrag{b1}{$b_1$}
%\psfrag{b2}{$b_2$}
%\psfrag{umu}{$u_\mu(x)$}
%\psfrag{s5}{$s_5$}
%\psfrag{p+}{$\psi_+(s_5)$}
%\begin{figure}
%\begin{center}
%\includegraphics[width=\textwidth]{cwfig1.eps}
%\caption[The Chacon-Walsh construction]{\label{fig:cwfig1}The above plot of $u_\mu$ shows some
%additional lines. Each line represents a step in the construction ---
%starting from zero, we run until we hit $a_1$ or $b_1$. If we hit $a_1$ first, we then run until we hit $a_2$ or $b_2$ The infimum of the old
%potential and the line gives our new potential. In the limit we aim to have potential agreeing with the target distribution.}
%\end{center}
%\end{figure}

We consider the embedding problem where we have a Brownian motion $B$
with $B_0 \sim \mu_0$ (an integrable starting distribution) and we
wish to embed an integrable target distribution $\mu$. This is
essentially the case considered by \citet{ChaconWalsh:76}, although
they only consider the case where $u_{\mu_0}(x) \ge u_\mu(x)$ for all
$x$ (when \eqref{eqn:pmulimit} implies $\mu_0$ and $\mu$ have the same
mean) --- we will see that this case is simpler than the general case
we consider. The embedding problem is frequently considered when
$\mu_0$ is the Dirac measure at $0$. One of the appealing properties
of the case where $B_0 = 0$ is that for all centred target
distributions \citep[Lemma~2.1]{Chacon:77}
\begin{equation}\label{eqn:modx}
u_\mu(x) \le -|x| = u_{\mu_0}(x),
\end{equation} 
and the condition on the ordering of potentials is easily satisfied.

We extend the technique of \cite{ChaconWalsh:76} to allow balayage on semi-infinite intervals. This extra step in the construction allows further flexibility later when we take limits of the constructions. In particular it will make the application of subsequent results trivial. Each step in the construction is described mathematically by a simple
balayage technique:
\begin{definition} \label{def:balayage}
Let $\mu$ be a probability measure on $\R$, and $I$ a finite, open
interval, $I=(a,b)$. Then define the {\em balayage} $\mu_I$ of $\mu$
on $I$ by:
\begin{eqnarray*}
\mu_I(A) & = & \mu(A) \hspace{2cm} A \cap \bar{I} = \emptyset;\\
\mu_I(\{a\}) & = & \int_{\bar{I}}\frac{b-x}{b-a}\, \mu(dx);\\
\mu_I(\{b\}) & = & \int_{\bar{I}}\frac{x-a}{b-a}\, \mu(dx);\\
\mu_I(I) & = & 0.
\end{eqnarray*}
Suppose now $I=(a,\infty)$ (resp.\ $I=(-\infty,a)$), and define the {\em
balayage} $\mu_I$ of $\mu$ by
\begin{eqnarray*}
\mu_I(A) & = & \mu(A) \hspace{2cm} A \cap \bar{I} = \emptyset;\\
\mu_I(\{a\}) & = & \int_{\bar{I}} \,\mu(dx);\\
\mu_I(I) & = & 0.
\end{eqnarray*}
\end{definition}

The balayage $\mu_I$ is a probability measure and if $I$ is a finite interval the means of $\mu$ and $\mu_I$ agree. In particular, $\mu_I$ is the
law of a Brownian motion started with distribution $\mu$ and run until
the first exit from $I$.

Our reason for introducing the Balayage technique is that the
potential of $\mu_I$ is readily calculated from the potential of
$\mu$:
\begin{lemma}[\citet{Chacon:77} Lemma~8.1] \label{lem:balayage}
Let $\mu$ be a probability measure with finite potential, $I=(a,b)$ a
finite open interval and $\mu_I$ the balayage of $\mu$ with respect to
$I$. Then
\begin{enumerate}
\item
$u_\mu(x) \ge u_{\mu_I}(x) \hspace{2cm} x \in \R$;
\item
$u_{\mu}(x) = u_{\mu_I}(x) \hspace{2cm} x \in I^C$;
\item
$u_{\mu_I}$ is linear for $x \in \bar{I}$.
\end{enumerate}
\end{lemma}

When $I$ is a semi-infinite interval we may calculate the potential in a similar way:

\begin{lemma} \label{lem:balayage2}
Let $\mu$ be a probability measure with finite potential $u_\mu$, $I =
(-\infty,a)$ or $I=(a,\infty)$ a semi-infinite interval and $\mu_I$
the balayage of $\mu$ with respect to $I$. Then
\begin{equation*}
\begin{array}{rcll}
u_{\mu_I}(x) & = & u_\mu(x) + \Delta m & x \notin I;\\
u_{\mu_I}(x) & = & u_\mu(a) + \Delta m - |a-x| & x \in I,
\end{array}
\end{equation*}
where we have written
\begin{equation*}
\Delta m = \int_I |x-a| \, \mu(dx).
\end{equation*}
\end{lemma}

The semi-infinite balayage step in Definition~\ref{def:balayage} can be recreated
using the balayage steps on compact intervals, for example
by taking the sequence of intervals
$(a,a+1),(a,a+2),(a,a+3),\ldots$. However it does not let us make the
same constructions as we can with the extended definition --- for example
if we wish our first step to be to move up to $1$, we would not be
able to carry out any further steps.

Formally, we may use balayage to define an embedding as the following
result shows. In the formulation of the result we assume we are given the sequence of functions we use to construct the stopping time, and from these deduce the target distribution. However we will typically use the result in situations where we have a desired target distribution and choose the sequence to fit this distribution.

\begin{lemma} \label{lem:fnembed2}
Let $f_1,f_2,\ldots$ be a sequence of linear functions on $\R$ such
that $|f_n'(x)|\le 1$ and
\begin{equation} \label{eqn:gefn}
g(x) = \inf_{n\in \N} f_n(x) \wedge (u_{\mu_0}(x)).
\end{equation}
Set $T_0=0$, $g_0(x)=u_{\mu_0}(x)$ and, for $n \ge 1$, define
\begin{eqnarray*}
a_n & = & \inf\{x \in \R: f_n(x) < g_{n-1}(x)\};\\
b_n & = & \sup\{x \in \R: f_n(x) < g_{n-1}(x)\};\\
T_n & = & \inf\{t \ge T_{n-1} : B_t \not\in (a_n,b_n)\};\\
g_n(x) & = & g_{n-1}(x) \wedge f_n(x).
\end{eqnarray*}
Then the $T_n$ are increasing so we define $T = \lim_{n \to \infty}
T_n$. If
\begin{equation} \label{eqn:gcond}
g(x) = u_\mu(x) -C
\end{equation}
for some $C \in \R$ and some integrable probability measure $\mu$ then
$T < \infty$ \as{} and $T$ is an embedding of $\mu$.
\end{lemma}

If we only consider the theorem under the condition $|f_n'(x)| < 1$ this is a formalised statement of the construction implicit in \citet{Chacon:77} and made explicit under the further condition \eqref{eqn:introcond} in \citet{ChaconWalsh:76}. The introduction of the balayage steps on the half-line is the novel content of the result.

\begin{proof}
The hard part is to show that if \eqref{eqn:gcond} holds then the
stopping time $T$ is almost surely finite. We prove in fact that
$\E(L_{T}) < \infty$, where $L$ is the local time of $B$ at zero. By
considering the martingale $|B_t| - L_t$ we must have
\begin{equation} \label{eqn:localtime}
\E(L_{T_1}) = u_{\mu_0}(0) - f_1(0) \wedge u_{\mu_0}(0)
\end{equation}
when the interval $(a_1,b_1)$ is compact; by approximating the semi-infinite interval by compact intervals (and a monotone convergence argument) this will extend to all possible choices of $f_1$, and, by an induction argument, we deduce
\[
\E(L_{T_n}) = u_{\mu_0}(0) - \inf_{k \le n} f_k(0) \wedge u_{\mu_0}(0).
\]
A monotone convergence argument allows us to deduce that
\[
\E(L_{T}) = u_{\mu_0}(0) - \inf_{n \in \N} f_n(0) \wedge u_{\mu_0}(0)
\]
which is finite by \eqref{eqn:gcond}, and hence $T < \infty$ \as{}.

The functions $g_n$ correspond to a potential of a measure $\mu_n$ ($\mu_n$ being the law of $B_{T_n}$) via:
\[
g_n(x) = u_{\mu_n}(x) - C_n
\]
for some constant $C_n$, and hence we have
\[
u_{\mu_0}(x) \ge u_{\mu_n}(x) - C_n \ge u_\mu(x) - C
\]
and as $n \to \infty$ the last two terms converge. From \eqref{eqn:pmulimit} we can deduce that 
\[
C_n \ge \left| \int x \, \mu_0(dx) - \int x \, \mu_n(dx)\right| \ge 0,
\]
so that since $g_n(x) \ge g(x)$, we have
\[
u_{\mu}(x) - C \le u_{\mu_n}(x) \le 0,
\]
the second inequality coming from the definition of the potential. Consequently, taking $x=0$, we can find a subsequence $n_j$ for which $\lim_{j \to \infty}u_{\mu_{n_j}}(0)$ exists, and hence for which $\lim_{n \to \infty} C_{n_j}$ also exists. Since $g_n(x)$ converges pointwise to $g(x)$ we must also have pointwise convergence of $u_{\mu_{n_j}}(x)$ to $u_{\mu}(x)-C'$ for some constant $C'$, and by Lemma~\ref{lem:limexists} $B_{T_{n_j}}$ converges weakly to $\mu$. Since also $T_{n_j} \uparrow T$, by the continuity of the Brownian motion $B_{T}$ has law $\mu$.

%Lemma~\ref{lem:balayage} implies
%\begin{equation*}
%u_{\mu_n}(x) =f_n(x) \wedge u_{\mu_{n-1}}(x) = \inf_{k\le n} f_k(x)
%\wedge (u_{\mu_0(x)}).
%\end{equation*}
%Since the functions $f_n$ satisfy
%\eqref{eqn:limitcond}, we know that the conditions of
%Lemma~\ref{lem:limexists} hold, determining the
%(unique) limiting distribution. Since $T_n \uparrow T< \infty$ \as{},
%$B_T$ has distribution $\mu$ --- \ie{} $T$ embeds $\mu$ --- by the
%continuity of the Brownian motion.
\end{proof}

The case considered by \citet{ChaconWalsh:76} has a notable
property. When the starting and target measures are centred (or at
least when their means agree) and
\begin{equation} \label{eqn:grcondition}
u_{\mu_0}(x) \ge u_\mu(x)
\end{equation}
then we may choose a construction such that $C=0$ in
\eqref{eqn:gcond}. In this case the process $B_{t \wedge T}$ is uniformly
integrable \citep[Lemma 5.1]{Chacon:77}. The desire to find a
condition to replace uniform integrability in situations where
\eqref{eqn:grcondition} does not hold, and to construct suitable
stopping times using this framework, is the motivation behind the
subsequent work.

We note also that --- for given $\mu, \mu_0$ --- we may find a construction for any $C$
which satisfies $C \ge \sup_{x} \left\{ u_{\mu}(x) - u_{\mu_0}(x)\right\}$; as we shall see, the case where there is equality is of particular interest. As
a consequence of \eqref{eqn:pmulimit} we must always have $C \ge 0$.

%We can now apply the graphical routine used before, along with the new
%`move' introduced of drawing the line to (plus or minus) infinity with
%gradient plus or minus $1$. This construction is shown in
%Figure~\ref{fig:cwfig5}.
%
%\begin{figure}[t]
%\begin{center}
%\includegraphics[width=\textwidth]{cwfig5.eps}
%\caption[The Chacon-Walsh construction for a shifted potential]{\label{fig:cwfig5}The above plot shows a potential $u_\mu$
%shifted so that it is no longer tangential to $-|x|$ at either
%$-\infty$ or $\infty$. This allows steps in the construction with
%gradients $\pm1$ as shown.}
%\end{center}
%\end{figure}
%
%The extended Chacon-Walsh embedding gives us a relatively large class
%of embeddings that (as a consequence of Remark~\ref{rem:ubound}) can
%be constructed for any integrable distributions $\mu, \mu_0$. We now
%turn to the question of which of these embeddings --- for given
%starting and target distributions --- are minimal.
%****************************************************************
%****************************************************************
%****************************************************************
%****************************************************************

%% file: minbasic.tex
\section{Minimality: Some Preliminary Results} \label{sec:minbasic}

In this and the subsequent section we discuss necessary and sufficient
conditions for an embedding of an integrable target distribution to be
minimal (Definition~\ref{def:minimal}) when we have an integrable starting distribution. These
results will extend the the conditions of Theorems~\ref{thm:monroeB}
and the following result:

\begin{theorem}[\citet{CoxHobson:03}]
\label{thm:mainB} Let $T$ be a stopping time of Brownian motion
which embeds an integrable distribution $\mu$ where $m = \int_\R x \, 
\mu(dx)
< 0$. Then the following conditions are equivalent:
\begin{enumerate}
\item $T$ is minimal for $\mu$;
\item for all stopping times $R \le S \le T$,
\begin{equation}\label{eqn:Scondineq0}
\E(B_S | \Fc_R) \le B_R \text{\quad \as{};}
\end{equation}
%\item for all stopping times $S \le T$,
%\begin{equation}\label{eqn:Scondineq1B}
%  \E (B_T | \F_S) \le B_S \text{\quad \as{};}
%\end{equation}
%\item for all $\gamma >0$
%\begin{equation*}
%\E(B_T;T>H_{-\gamma}) \le -\gamma \Pr(T>H_{-\gamma});
%\end{equation*}
%\item as $\gamma \to \infty$
%\begin{equation*}
%\gamma \Pr(T > H_{-\gamma}) \to 0;
%\end{equation*}
%\item
%the family $\{B_S^-\}$ taken over stopping times $S \le T$ is uniformly 
%integrable;
%\item
%for all $x > 0$
%\begin{equation*}
%\E(B_{T \wedge H_x}) = 0.
%\end{equation*}
\end{enumerate}
In the case where $\supp(\mu) \subseteq [\alpha,\infty)$ for some
$\alpha < 0$ then the above conditions are also equivalent to the
condition:

\noindent
{\it(iii)}
\begin{equation}\label{eqn:TleqHaA}
  \Pr(T \le H_\alpha) = 1,
\end{equation}
where $H_\alpha = \inf \{ t \ge 0 : B_t = \alpha\}$ is the hitting time of $\alpha$.
\end{theorem}

\begin{remark} For further necessary and sufficient conditions, see also \citet{CoxHobson:03}.
\end{remark}

As a starting point, we recall:
\begin{prop}[\citet{Monroe:72}, Proposition 2] \label{prop:existenceB} For 
any stopping time $T$ there exists a minimal stopping time $S 
\le T$ 
such that $B_S \sim B_T$. \end{prop}

Monroe's proof does not rely on the fact that $B$
starts at $0$, and so the result extends to a general starting
distribution.

% the proof of
%Proposition~\ref{prop:existenceB} does not rely on the fact that $B$
%starts at $0$, and so the result extends to a general starting
%distribution, so that there always exists a minimal embedding smaller
%than any given embedding.

It can also be seen that the argument used in \citet{Monroe:72} to
show that if the process is uniformly integrable then the process is
minimal does not require the starting measure to be a point mass. For
completeness we state a similar result, with the proof given in \citet{Monroe:72}:
\begin{lemma}
Let $T$ be a stopping time embedding $\mu$ in $(B_t)_{t \ge 0}$, with
$B_0 \sim \mu_0$ where $\mu$ and $\mu_0$ are integrable distributions. If
\begin{equation}\label{eqn:monroenew}
\E(B_T|\Fc_S) = B_S \mbox{ \as{}}
\end{equation}
for all stopping times $S \le T$ then $T$ is minimal.
\end{lemma}
Note that $S \equiv 0$ implies that $\mu, \mu_0$ have the same mean.
%\begin{proof}
%Let $S \le T$ be a stopping time such that $B_S \sim \mu$. Then for $a
%\in \R$
%\begin{equation*}
%\E(B_T ; B_T \ge a) = \E(B_S ; B_S \ge a) = \E(B_T ; B_S \ge a).
%\end{equation*}
%Consequently $B_S = B_T$ \as{}. If $R$ is another stopping time, $S
%\le R \le T$, then
%\begin{equation*}
%B_R = \E(B_T | \Fc_R) = \E(B_S| \Fc_R) = B_S = B_T \mbox{ \as{}}.
%\end{equation*}
%And by the continuity of Brownian paths, $B$ is constant on the
%interval $[S,T]$ and hence $S=T$ \as{}.
%\end{proof}

\begin{remark} \label{rem:simpleexample}
We will later be interested also in necessary conditions for
minimality. The condition in \eqref{eqn:monroenew} is not necessary
even when both starting and target measures are centred, as can be
seen by taking $\mu_0 = \half \delta_{-1} + \half \delta_{1}$ and $\mu
= \delta_0$, where it is impossible to satisfy \eqref{eqn:monroenew}
but the (only) minimal stopping time is `stop when the process hits
0.'
\end{remark}

The condition in \eqref{eqn:monroenew} is equivalent to uniform
integrability of the process $(B_{t \wedge T})_{t \ge 0}$. One
direction follows from the optional stopping theorem, the reverse
implication comes from the upward martingale theorem
\citep{RogersWilliams:00a}[Theorem~II.69.5], which tells us that the
process $X_t = \E(B_T|\Fc_t)$ is a uniformly integrable martingale on
$t \le T$. When \eqref{eqn:monroenew} holds, $X_t = B_{t \wedge T}$,
and the process $B_{t \wedge T}$ is a uniformly integrable martingale.

For the rest of this section we will consider minimality for general
starting and target measures: particularly when the means do not
agree. If this occurs when the starting measure is a point mass,
necessary and sufficient conditions are given in
Theorem~\ref{thm:mainB}. In subsequent proofs with general starting
measures we will often reduce problems to the point mass case in order
to apply the result.

\begin{remark} \label{rem:bounded}
The condition given in {\it(iii)} of Theorem~\ref{thm:mainB} hints at a
more general idea inherent in the study of embeddings in Brownian
motion. When $B_0 = 0$, it is a well known fact that if there exists
$\alpha < 0 < \beta$ such that $T \le H_\alpha \wedge H_\beta$ then
$B_{t \wedge T}$ is a uniformly integrable martingale. If $T \le
H_\alpha$ then the process is a supermartingale. In terms of
embeddings, this observation has the following consequence: if the
target distribution is centred and supported on a bounded interval, an
embedding is minimal if and only if the process never leaves this
interval. If the target distribution has a negative mean, but still
lies on a bounded interval, any embedding must move above the interval
--- \ie{} $\Pr(\sup_{t \le T} B_t \ge x) > 0$ for all $x \ge
0$. Theorem~\ref{thm:mainB} and Proposition~\ref{prop:existenceB} tell
us that in this case an embedding exists for which $T \le H_\alpha$
and all minimal embeddings satisfy this property.
\end{remark}

Recall that there is a natural ordering on the set of (finite)
measures on $\R$, that is $\mu \preceq \nu$ if and only if $\mu(A) \le
\nu(A)$ for all $A \in \Borel(\R)$, in which case we say that $\nu$
dominates $\mu$. In such instances it is possible to define a
(positive, finite) measure $(\nu-\mu)(A) = \nu(A) - \mu(A)$. The
notation $\nu = \Lc(B_T;T < H_\alpha)$ is used to mean the
(sub-probability) measure $\nu$ such that $\nu(A) = \Pr( B_T \in A, T
< H_\alpha)$.

\begin{lemma} \label{lem:reduceminimal}
Let $B_t$ be a Brownian motion with $B_0=0$, $T$ a stopping time
embedding a distribution $\mu$, $\tilde{\mu}$ a target distribution
such that $\supp(\tilde{\mu}) \subseteq [\alpha,\infty)$ for some
$\alpha < 0$ and $\int x \,\tilde{\mu}(dx) \le 0$. Then if $\nu =
\Lc(B_T;T < H_\alpha)$ is dominated by $\tilde{\mu}$, there exists a
minimal stopping time $\tilde{T} \le T \wedge H_\alpha$ which embeds
$\tilde{\mu}$.

Similarly, if $\tilde{\mu}$ is such that $\supp(\tilde{\mu}) \subseteq
[\alpha,\beta]$ and $\int x \, \tilde{\mu}(dx) = 0$, and if
$\nu=\Lc(B_T ; T < H_{\alpha} \wedge H_{\beta})$ is dominated by
$\tilde{\mu}$, then there exists a minimal stopping time $\tilde{T}
\le T \wedge H_{\alpha} \wedge H_{\beta}$ which embeds $\tilde{\mu}$.
\end{lemma}

\begin{proof}
Construct a stopping time $T'$ as follows: on $\{ T < H_\alpha\}$,
$T'=T$; otherwise choose $T'$ so that $T' = H_\alpha + T'' \circ
\theta_{H_\alpha}$ where $T''$ is chosen to embed $(\tilde{\mu} -
\nu)$ on $\{T' \ge H_\alpha\}$ given $B_0 = \alpha$. Then $T'$ is an
embedding of $\tilde{\mu}$ and $T' \le T$ on $\{T < H_\alpha\}$. So by
Proposition~\ref{prop:existenceB} we may find a minimal embedding
$\tilde{T} \le T' \wedge H_\alpha = T \wedge H_\alpha$ which embeds
$\tilde{\mu}$.

The proof in the centred case is essentially identical, but now
stopping the first time the process leaves $[\alpha, \beta]$.
\end{proof}

We turn now to the case of interest --- that is when $B_0 \sim \mu_0$
and $B_T \sim \mu$ for integrable measures $\mu_0$ and $\mu$. The
following lemma is essentially technical in nature, but will allow us
to deduce the required behaviour on letting $A$ increase in density.

\begin{lemma} \label{lem:Acond}
Let $T$ be a minimal stopping time, and $A$ a countable subset of $\R$
such that $A$ has finitely many elements in every compact subset of
$\R$ and $d(x,A) < M$ for all $x \in \R$ and some $M>0$. We consider
the stopping time
\begin{equation*}
R(A) = \inf \{ t \ge 0 : B_t \in A \} \wedge T
\end{equation*}
and we write
\begin{equation*}
E_A(x) = 
\begin{cases}
\E(B_T| T > R(A), B_{R(A)} = x)&: \Pr(T > R(A), B_{R(A)} = x) >0;\\
x&: \Pr(T > R(A), B_{R(A)} = x) =0.
\end{cases}
\end{equation*}
Then there exists $a \in \bar{\R} = \R \cup \{-\infty\} \cup \{
\infty\}$ such that
\begin{eqnarray}
E_A(x) > x & \implies & x < a, \label{eqn:EAless}\\
E_A(x) < x & \implies & x >a, \label{eqn:EAmore}
\end{eqnarray}
and $T \le H_a$ on $\{T\ge R(A)\}$.

Further, if there exists $x<y$ such that $E_A(x) > x$ and $E_A(y)<y$
then there exists $a_\infty \in [x,y]$ such that $T\le H_{a_\infty}$.
\end{lemma}

\begin{proof}
Suppose that there exists $x<y$ such that $E_A(x) < x$ and $E_A(y) >
y$, and suppose $E_A(w) = w$ for $x<w<y$. We show that we can
construct a strictly smaller embedding, contradicting the assumption
that $T$ is minimal.

Define the stopping time $T' = R(A) \indic{B_{R(A)} \in \{x,y\}} +
T\indic{B_{R(A)} \notin \{x,y\}}$ and for some $z\in (x,y)$, the
stopping time
\begin{equation*}
T'' = \inf\{ t \ge T': B_t = z \} \wedge T.
\end{equation*}
As a consequence of Remark~\ref{rem:bounded}, paths from both $x$ and
$y$ must hit $z$.

Consider the set $\{T'' < T\}$. On this set we have only paths with
$B_{R(A)} = x$ and $B_{R(A)} = y$. Define $\mu_x = \Lc(B_T ; B_{R(A)}
= x, T''<T)$ and $\mu_y = \Lc(B_T; B_{R(A)} = y, T''<T)$. Since
Brownian motion bounded above is a submartingale,
\begin{equation*}
\E(B_{T \wedge H_z}; B_{R(A)} = x, T>R(A)) \ge x \Pr(B_{R(A)}=x, T >R(A)).
\end{equation*}
Together with $E_A(x) < x$ this implies
\begin{equation*}
z \Pr (B_{R(A)} = x, T'' < T) > \E(B_T ; B_{R(A)} = x, T'' < T),
\end{equation*}
\ie{} we must have $\frac{1}{\mu_x(\R)}\int w \, \mu_x(dw)<z$, and
similarly $\frac{1}{\mu_y(\R)}\int w \, \mu_y(dw)>z$. Then we apply
Lemma~\ref{lem:reduceminimal} to the processes $B_{T''+t}$ on
$\{B_{R(A)} = x, T''<T\}$ and $\{B_{R(A)} = y, T''<T\}$ with the
measures
\begin{eqnarray*}
\tilde{\mu}_x & = & \mu_x|_{[a_1\infty)} + \mu_y|_{(a_2,\infty)}\\
\tilde{\mu}_y & = & \mu_x|_{(-\infty,a_1)} + \mu_y|_{(-\infty,a_2)}
\end{eqnarray*}
where we choose $a_1<z<a_2$ so that 
\begin{equation*}
\frac{1}{\tilde{\mu}_x(\R)} \int
w \, \tilde{\mu}_x(dw) \le z \mbox{ and } \frac{1}{\tilde{\mu}_y(\R)} \int w
\, \tilde{\mu}_y(dw) \ge z
\end{equation*}
and also so that $\mu_x(\R) = \tilde{\mu}_x(\R)$ and $\mu_y(\R) =
\tilde{\mu}_y(\R)$\footnote{It may be necessary to consider only a
proportion of the paths hitting $z$ from one side; this can be done by
choosing paths according to an independent $U([0,1])$ random variable
and running the rest of the paths according to $T$. This will still
construct a strictly smaller stopping time.}. This will produce a
strictly smaller embedding, in contradiction to the assumption that
$T$ is minimal.

So we have shown that there exists $a$ such that \eqref{eqn:EAless}
and \eqref{eqn:EAmore} hold. We just need to show that we can choose
$a$ so that $T \le H_a$ on $\{T \ge R(A)\}$.

Suppose that there exists $x<y$ such that $E_A(x)>x$ and $E_A(y)<y$
and $E_A(w) = w$ for $w \in (x,y)$. If
\begin{equation} \label{eqn:suppcond}
\sup_{x < a} \mu_x((a,\infty))
= 0 \mbox{ and } \sup_{y > a} \mu_y((-\infty,a)) = 0
\mbox{ for some } a \in (x,y)
\end{equation}
then $T$ minimal and Theorem~\ref{thm:mainB} implies that $T \le H_a$
on $\{T \ge R(A)\}$.

So suppose that \eqref{eqn:suppcond} does not hold. We shall show that
we can then find a sequence $x_1, x_2, \ldots, x_r$ of elements of $A$ such
that we are able to transfer mass between the $x_i$ to produce a
smaller embedding. We begin by choosing $x_1$ to be the point of $A$
satisfying $E_A(x)>x$ for which the support of $\mu_{x}$ extends
furthest to the right, and $y_1$ similarly the point satisfying
$E_A(y)<y$ for which the support of $\mu_{y}$ extends furthest to the
left. If the support of these measures overlap we show we can exchange
mass between $\mu_{x_1}$ and $\mu_{y_1}$ and embed to find a smaller
stopping time. Otherwise we look at those points for which $E_A(x) =
x$ and the support overlaps that of $\mu_{x_1}$ but extends further to
the right. In this way we can find a sequence whose supports overlap
(since \eqref{eqn:suppcond} does not hold) and we may again perform a
suitable exchange of mass to show that we can find a smaller
embedding. Then we take $x_r = y_1$ and the points satisfy $x_2 < x_3
< \ldots < x_{r-1}$.

There are several technical issues we need to address. Firstly, if we
find at some stage there are two points which both satisfy the
criterion --- for example their supports have the same upper bound ---
then we may use either point. Secondly, if the support of all suitable
points has a maximum which is not attained we may still use the same
procedure but we must (and can) choose a point which approximates the
bound suitably closely for subsequent steps to work. Finally we note
that once we choose $x_2$, since there is at most one point to the
right of $y_1$, there exists only a finite number of points left to
choose from (by assumption on $A$) and so the sequence will be finite.

The technical construction is as follows: let $x_1$ be the largest
value such that $E_A(x_1)>x_1$ and
\begin{equation*}
\sup\{z:z \in \supp(\mu_{x_1})\} =
    \sup_{w: E_A(w)>w} \{\sup\{z:z \in \supp(\mu_{w})\}\},
\end{equation*}
(or at least so that the left hand side approximates the right hand
side sufficiently closely for the next step to work --- since the support of the points to the right overlaps we shall be able to find $x_1$ with supremum of its support sufficiently close to the term on the left) and let $y_1$ be
the smallest value such that $E_A(y_1)<y_1$ and
\begin{equation*}
\inf\{z:z \in \supp(\mu_{y_1})\} =
    \inf_{w: E_A(w)<w} \{\inf\{z:z \in \supp(\mu_{w})\}\}.
\end{equation*}
Then (by the assumption that \eqref{eqn:suppcond} does not hold) we
can find a sequence $x_1,x_2,\ldots,x_r$ such that $x_r =y_1$ and $x_2
< x_3 < \ldots x_{r-1}$, $E_A(x_i) = x_i$ for $1 < i < r$ and, if we
define $I_i = \inf \{ \mbox{intervals } I : \supp(\mu_{x_i}) \subseteq
I\}$, then
\begin{eqnarray}
Leb(I_i \cap I_{i+1}) &>& 0 \hspace{2cm} k=1,\ldots,r-1,\nonumber\\
Leb(I_i \cap I_{i+2}) &=& 0 \hspace{2cm} k=1,\ldots,r-2.\label{eqn:nointersect}
\end{eqnarray}
This is done by choosing at each step the $w$ with $E_A(w) = w$ which
overlaps the support of the previous $\mu_{x_i}$ and whose support
extends furthest to the right, until the support overlaps with the
support of $\mu_{y_1}$.

We write $\mu_i = \mu_{x_i}$. For general $1 \le i < r$ now consider
$\mu_i'$ defined by
\begin{equation*}
\mu_i' = \mu_i |_{(-\infty,y_i)} + \mu_{i+1}|_{(-\infty,y_i)}
\end{equation*}
where $y_i$ is chosen such that $\mu_i([y_i,\infty)) =
\mu_{i+1}((-\infty,y_i))$. Then it must be true that $\int w \,
\mu_i'(dw) < \int w \, \mu_i(dw)$. Define
\begin{eqnarray*}
m_i & = & \int w \, \mu_i(dw) - \int w \, \mu_i'(dw) >0\\
m_0 & = & \int w \, \mu_1(dw) - \mu_1(\R) x_1 > 0\\
m_r & = & \mu_r(\R) x_r - \int w \, \mu_r(dw) > 0 
\end{eqnarray*}
and set $\Delta m = \inf \{ m_i : 0 \le i \le r\}$. Then for each $i$
we can find $v_i < z_i$ such that $\mu_i([z_i,\infty)) =
\mu_{i+1}((-\infty,v_i))$ and for
\begin{equation*}
\mu_i' = \mu_i |_{(-\infty,z_i)} + \mu_{i+1}|_{(-\infty,v_i)}
\end{equation*}
we have
\begin{equation*}
\int w \,\mu_i(dw) - \int w \,\mu_i'(dw) = \Delta m.
\end{equation*}
Set
\begin{eqnarray*}
\mu_1'' & = & \mu_1 |_{(-\infty,z_1)} + \mu_2 |_{(-\infty,v_1)},\\
\mu_i'' & = & \mu_{i-1}|_{[z_{i-1},\infty)} + \mu_{i} |_{[v_{i-1},z_i)}
      + \mu_{i+1} |_{(-\infty,v_i)} \hspace{2cm} i=2,\ldots,r-1,\\
\mu_r'' & = & \mu_{r-1} |_{[z_{r-1},\infty))}
             + \mu_r |_{[v_{r-1},\infty)}.
\end{eqnarray*}
Then
\begin{eqnarray*}
\int x \, \mu_1'' & \ge & \mu_1''(\R)x_1\\
\int x \, \mu_i'' & = & \mu_i''(\R)x_i \hspace{2cm} i=2,\ldots,r-1\\
\int x \, \mu_r'' & \le & \mu_r''(\R)x_r.
\end{eqnarray*}
So the conditions of Lemma~\ref{lem:reduceminimal} are satisfied (due to \eqref{eqn:nointersect}) for
each $\mu_i''$ and we can find strictly smaller stopping times on each
of the sets $\{T>R(A), R(A) = x_i\}$.

It only remains to show the final statement of the lemma.  Let $A'
\supset A$ be another set satisfying the conditions of the lemma for
some $M'$, such that $A'\setminus A \subseteq [x,y]$. Then there
exists $x',y' \in A'$ such that $x \le x'<y' \le y$, $E_{A'}(x') > x'$
and $E_{A'}(y') < y'$ --- if this were not the case at least one of
the embeddings conditional on $\{R(A') = z\}$ would not be minimal.

Now consider a sequence $A \subset A_1 \subset A_2 \subset \ldots$ and
such that $A_n \setminus A \subseteq [x,y]$ and $d(z,A_n) \le 2^{-n}$
for $z \in [x,y]$. Let
\begin{eqnarray*}
\Lambda & = & \{ a \in [x,y] : T \le H_a \mbox{ on } \{ T \ge R(A) \} \};\\
\Lambda_n & = & \{ a \in [x,y] : T \le H_a \mbox{ on } \{ T \ge R(A_n) \} \}.
\end{eqnarray*}
Then the sets $\Lambda,\Lambda_n$ are closed, $\Lambda \supseteq
\Lambda_1 \supseteq \Lambda_2 \supseteq \ldots$, and each $\Lambda_n$
is non-empty. So there exists $a_\infty \in \Lambda_n$ for all
$n$. Hence $T \le H_{a_\infty}$ on $\{T \ge R(A_n)\}$ for all $n$. But
$R(A_n) \downarrow 0$ on $\{B_0 \in [x,y]\}$ and $R(A) \le
H_{a_\infty}$ on $\{B_0 \not\in [x,y]\}$.
\end{proof}

This result, although technical in nature, can be thought of as
beginning to describe the sort of behaviour we shall expect from
minimal embeddings in this general context. The cases considered in
Theorem~\ref{thm:mainB} suggest behaviour of the form: `the process
always drifts in the same direction', if indeed it drifts at all. The
example of Remark~\ref{rem:simpleexample} suggests that this is not
always possible in the general case, and the previous result suggests
that this is modified by breaking the space into two sections, in each
of which the process can be viewed separately. The way these sections
are determined is clearly dependent on the starting and target
measures, and we shall see in the next section that the potential of
these measures provides an important tool in determining how this
occurs.

%****************************************************************
%****************************************************************
%****************************************************************
%****************************************************************

%% file: minandpot.tex
\section{Minimality and Potential} \label{sec:minandpot}

The main aim of this section is to find equivalent conditions to
minimality which allow us to characterise minimality simply in terms
of properties of the process $B_{t \wedge T}$. This is partly in order
to prove the following result:
\begin{quote}
The Chacon-Walsh type embedding is minimal when constructed using the
functions $u_{\mu_0}$ and $c(x)=u_\mu(x) -C$ where 
\begin{equation} \label{eqn:uccond}
C = \sup_{x} \{u_\mu(x) - u_{\mu_0}(x)\}.
\end{equation}
\end{quote}

We have already shown that provided the means of our starting and
target distribution match, and \eqref{eqn:grcondition} holds (so that
$C=0$ --- the solution in this case to \eqref{eqn:uccond}), then the
process constructed using the Chacon-Walsh technique is uniformly
integrable, and therefore minimal. Of course the Chacon-Walsh
construction is simply an example of an embedding, and the functions
$u_{\mu_0}$ and $c$ are properties solely of the general problem ---
it seems reasonable however that these functions will appear in the general
problem of classifying all minimal embeddings.

So consider a pair $\mu_0,\mu$ of integrable
measures. Remark~\ref{rem:potprop}(iv) tells us we we can choose $C< \infty$ such
that \eqref{eqn:uccond} holds. We know $u_{\mu_0}(x)-c(x)$ is bounded
above, and $\inf_{x \in \R} u_{\mu_0}(x)-c(x) = 0$. We consider
\begin{equation} \label{eqn:scAdefn}
\scA = \{ x \in [-\infty,\infty] : \lim_{y \to x}[u_{\mu_0}(y)-c(y)] = 0\}.
\end{equation}
Since both functions are Lebesgue almost-everywhere differentiable,
Remark~\ref{rem:potprop}(v) implies $\scA \subseteq \scA'$ where
$\scA'$ is the set
\begin{equation} \label{eqn:Asubset}
\{ x \in [-\infty,\infty]: \mu((-\infty,x)) \le \mu_0((-\infty,x)) \le
\mu_0((-\infty,x]) \le \mu((-\infty,x])\}.
\end{equation}
One consequence of this is that if the starting distribution has an
atom at a point of $\scA$ then the target distribution has an atom at
least as large. Also we introduce the following definition. Given a
measure $\nu$, $a \in \R$ and $\theta \in
[\nu((-\infty,a)),\nu((-\infty,a])]$ we define the measure
$\check{\nu}^{a,\theta}$ to be the measure which is $\nu$ on
$(-\infty,a)$, has support on $(-\infty,a]$ and $\check
\nu^{a,\theta}(\R)=\theta$. We also define $\hat{\nu}^{a,\theta} =
\nu-\check{\nu}^{a,\theta}$. Then for $a \in \scA$ we may find
$\theta$ such that
\begin{eqnarray*}
\check{\mu}^{a,\theta}((-\infty,a]) & = & \check{\mu}^{a,\theta}_0((-\infty,a])\\
\hat{\mu}^{a,\theta}([a,\infty)) & = & \hat{\mu}^{a,\theta}_0([a,\infty)).
\end{eqnarray*}
When $\mu_0((-\infty,a)) < \mu_0((-\infty,a])$ there will exist
multiple $\theta$. We will occasionally drop the $\theta$ from the
notation since this is often unnecessary.

These definitions allows us to write the potential in terms of the new
measures (for any suitable $\theta$)
\begin{equation} \label{eqn:umuis}
u_{\mu}(x) = \int_{(-\infty,x]}(y-x)\, \check{\mu}^x(dy) + \int_{[x,\infty)} (x-y) \hat{\mu}^x(dy).
\end{equation}
As a consequence of this and a similar relation for $u_{\mu_0}$, we
are able to deduce the following important facts about the set $\scA$:
\begin{itemize}
\item
if $x<z$ are both elements of $\scA$ (possibly $\pm\infty$), then
\begin{equation} \label{eqn:meansequal}
\int y \, (\mu - \check{\mu}^{x,\theta} - \hat{\mu}^{z,\phi})(dy) = \int y \, (\mu_0 - \check{\mu}^{x,\theta}_0 - \hat{\mu}_0^{z,\phi})(dy).
\end{equation}
That is, we may find measures agreeing with $\mu$ and $\mu_0$ on
$(x,z)$ and with support on $[x,z]$ which have the same mean.

\item If $x \in \scA$, by definition
\begin{equation} \label{eqn:uineq}
u_\mu(x) - u_{\mu_0}(x) \ge \lim_{z \to -\infty} (u_{\mu}(z) - u_{\mu_0}(z)).
\end{equation}
This can be rearranged, using \eqref{eqn:umuis}, to deduce
\begin{equation*}
\int_{(-\infty,x]} y \, \check{\mu}^x_0 (dy) \le \int_{(-\infty,x]} y \, \check{\mu}^x(dy)
\end{equation*}
with equality if and only if there is also equality in
\eqref{eqn:uineq} --- that is when $-\infty \in \scA$.
\end{itemize}

Together these imply that the set $\scA$ divides $\R$ into intervals
on which the starting and target measures place the same amount of
mass. Further, the means of the distributions agree on these intervals
except for the first (resp. last) interval where the mean of the
target distribution will be larger (resp. smaller) than that of the
starting distribution unless $-\infty$ (resp. $\infty$) is in $\scA$,
when again they will agree. Note the connection between this idea and
Lemma~\ref{lem:Acond}

Before we prove the result we establish several results that are
needed in the proof.

\begin{prop} \label{prop:cinfinA}
Suppose $T \le H_{a_\infty}$ is an embedding of $\mu$ for $a_\infty
\in \R$. Then $a_\infty \in \scA$.
\end{prop}
\begin{proof}
Clearly $a_\infty$ must lie in $\scA'$ (see \eqref{eqn:Asubset}).
Suppose also that $a_\infty< z \in \scA$. We may choose $\theta, \phi$
such that $\mu_0 - \check{\mu}^{a_\infty,\theta}_0 -
\hat{\mu}^{z,\phi}_0$ has no atom at either $a_\infty$ or $z$.

Then 
\begin{equation} \label{eqn:umuineq}
u_{\mu_0}(a_\infty) \ge u_{\mu}(a_\infty)-C
\end{equation}
and $C = u_{\mu}(z) - u_{\mu_0}(z)$ imply
\begin{equation} \label{eqn:muhatineq}
\int y \, (\mu - \check{\mu}^{a_\infty,\theta} - \hat{\mu}^{z,\phi})(dy) \ge \int y \, (\mu_0 - \check{\mu}^{a_\infty,\theta}_0 - \hat{\mu}^{z,\phi}_0)(dy),
\end{equation}
the term on the right being equal to $\E(B_0 ; B_0 \in (a_\infty,z))$
and the term on the left at most $\E(B_T ; B_0 \in
(a_\infty,z))$. However $B_T = B_{T \wedge H_{a_\infty}}$ is a
supermartingale on $\{B_0 \ge a_\infty\}$, so we must have equality in
\eqref{eqn:muhatineq} and hence in \eqref{eqn:umuineq}. So $a_\infty
\in \scA$. A similar proof can be used for $z<a_\infty$. If there does not exists any such $z$, $a_\infty \in \scA$, since $\scA \neq \emptyset$.
\end{proof}

\begin{prop} \label{prop:EAgx}
Suppose $T$ is minimal and $A$ is a countable subset of $\R$ such that
$A$ has finitely many elements in every compact subset of $\R$ and
$d(x,A) < M$ for all $x \in \R$ and some $M>0$. Suppose also that $S
\le T$ is a stopping time and $I \subseteq \R$ is an
interval such that $\partial I \subseteq A$. If
\begin{equation} \label{eqn:BTgBS}
\E(B_T ; F \cap \{B_0 \in I\}) > \E(B_S ; F \cap \{B_0 \in I\})
\end{equation}
for some $F \in \Fc_S$ then $E_A(x) > x$ for some $x \in A \cap \bar{I}$.
\end{prop}
\begin{proof}
We may assume $F \subseteq \{B_0 \in I\}$ and we note that therefore $B_{R(A)}
\in A \cap \bar{I}$ on $\{R(A) < T\} \cap F$. Since $B_{t \wedge R(A)}$ is uniformly
integrable,
\begin{eqnarray*}
\E(B_S ; F) & = & \E(B_{R(A)};F \cap \{S \le R(A)\}) + \E(B_S ; F \cap \{R(A) < S\})\\
\E(B_T ; F) & = & \E(B_{R(A)};F \cap \{T =  R(A)\}) + \E(B_T ; F \cap \{R(A) < T\}).
\end{eqnarray*}
So \eqref{eqn:BTgBS} and the above identities imply
\begin{equation*}
\begin{split}
\E(B_T;F \cap \{R(A)<T\}) >& \E(B_{R(A)} ; F \cap \{ S \le R(A) < T\}) \\&{}+ \E(B_S; F \cap \{R(A) < S\}).
\end{split}
\end{equation*}
However if $E_A(x) \le x$ for all $x \in A \cap \bar{I}$ and $T$ is
minimal, by Theorem~\ref{thm:mainB}:
\begin{eqnarray*}
\E(B_T; F \cap \{R(A) < S\}) & \le & \E(B_S ; F \cap \{R(A) < S\})\\
\E(B_T; F \cap \{S \le R(A) < T\}) & \le & \E(B_{R(A)} ; F \cap \{S \le R(A) < T\})
\end{eqnarray*}
and we deduce a contradiction.
\end{proof}

\begin{prop} \label{prop:impliesmg}
Suppose $F \in \Fc_0$, $\E(B_T;F) = \E(B_0;F)$ and
\begin{equation} \label{eqn:equalmg}
\E(B_T|\Fc_S) \le B_S \mbox{ on } F
\end{equation}
for all stopping times $S$. Then in fact we have equality --- that is
\begin{equation*}
\E(B_T|\Fc_S) = B_S
\end{equation*}
almost surely on $F$.
\end{prop}
\begin{proof}
If $\Pr(F)=0$ there is nothing to prove. Otherwise we may condition on
$F$ to reduce to showing the result when $F = \Omega$.

By the upward martingale theorem
\citep{RogersWilliams:00a}[Theorem~II.69.5], the process
\begin{equation*}
X_t = \E(B_T | \Fc_t)
\end{equation*}
is uniformly integrable. Also $\E(B_T | \Fc_0) \le B_0$ and $\E(B_T) =
\E(B_0)$ implies $\E(B_T|\Fc_0) = B_0$. Let $Y_t = B_{T \wedge t} -
X_{T \wedge t}$. By \eqref{eqn:equalmg} $Y_t$ is a non-negative local
martingale such that $Y_0 = Y_T =0$. Hence $Y \equiv 0$.
\end{proof}

\begin{lemma} \label{lem:genmin}
If $T$ is minimal and $a \in \scA$ then $T \le H_a$ and
\begin{eqnarray}
\E(B_T | \Fc_S) & \le & B_S \mbox{ on } \{B_0 \ge a\}; \label{eqn:Bdown}\\
\E(B_T | \Fc_S) & \ge & B_S \mbox{ on } \{B_0 \le a\}. \label{eqn:Bup}
\end{eqnarray}
\end{lemma}

\begin{proof}%[Proof of Lemma~\ref{lem:genmin}]
Suppose initially $a \in \R$. Let $\theta = \mu_0((-\infty,a))$. If
$\{B_0 < a\} \not\subseteq \{B_T \le a\}$ \as{} then also $\{B_0 \ge
a\} \not\subseteq \{B_T \ge a\}$ \as{} and
\begin{eqnarray*}
\E(B_0; B_0 < a) = \int y \, \check \mu_0^{a,\theta}(dy) & \le & \int y \, \check \mu^{a,\theta}(dy) < \E(B_T;B_0 < a);\\
\E(B_0; B_0 \ge a) = \int y \, \hat \mu_0^{a,\theta}(dy) & \ge & \int y \, \check \mu^{a,\theta}(dy) > \E(B_T;B_0 < a).
\end{eqnarray*}
So there exists $x_1 \le a$ and $x_2 \ge a$ such that (by
Proposition~\ref{prop:EAgx})
\begin{equation*}
E_A(x_1) < x_1 \mbox{ and } E_A(x_2) > x_2
\end{equation*}
for a suitable choice of $A$ --- a contradiction to
Lemma~\ref{lem:Acond}.

A similar argument can be used with $\theta = \mu_0((-\infty,a])$ to
deduce that $\{B_0 \le a\} \subseteq \{B_T \le a\}$ \as{} and $\{B_0
\ge a\} \subseteq \{B_T \ge a\}$ \as{}. So if there is an atom of
$\mu_0$ at $a$ then paths starting at $a$ must also stop at $a$, and
hence (by the minimality of $T$) must stop immediately --- \ie{} $T=0$
on $\{B_0 = a\}$.

So consider paths for which $\{B_0 < a\}$. For almost all these paths,
for some choice of $A$, $B_{R(A)} < a$. If \eqref{eqn:Bup} fails, by
Proposition~\ref{prop:EAgx} there exists $x < a$ such that $E_A(x) <
x$. Then Lemma~\ref{lem:Acond} and (for $\theta = \mu_0((-\infty,a))$)
\begin{equation*}
\int y \, \check \mu^{a,\theta}_0(dy) \le \int y \, \check\mu^{a,\theta}(dy)
\end{equation*}
imply there must also exist $y<x$ such that $E_A(y) > y$, and hence
$a'<a$ such that $T \le H_{a'}$. Then $B_{t \wedge T}$ is a
supermartingale on $\{B_0 > a'\}$ (and a submartingale on $\{B_0 \le
a'\}$). But Proposition~\ref{prop:cinfinA} and \eqref{eqn:meansequal}
imply $\E(B_0; a' < B_0 < a) = \E(B_T; a' < B_0 < a)$ and therefore
(by Proposition~\ref{prop:impliesmg}) $B_{t \wedge T}$ is a true
martingale on $\{a' < B_0 < a\}$ --- in particular $T \le H_a$ on
$\{B_0 < a\}$, and \eqref{eqn:Bup} holds. Similarly \eqref{eqn:Bdown}
can be shown to hold.

So suppose now that $a = \infty$ (the case $a=-\infty$ is similar) and
there exists $a' < \infty$ also in $\scA$. By the above, $T \le
H_{a'}$ and so $B_{t \wedge T}$ is a supermartingale on $\{B_0 >
a'\}$, while by \eqref{eqn:meansequal} $\E(B_0 ; B_0>a') = \E(B_T ;
B_0 > a')$, and hence $B_{t \wedge T}$ satisfies \eqref{eqn:Bup} by
Proposition~\ref{prop:impliesmg}.

Finally suppose $\scA = \{\infty\}$. By Lemma~\ref{lem:Acond} $E_A(x)
\ge x$ for all suitable choices of $A$ and all $x$. Hence, by
Proposition~\ref{prop:EAgx},
\begin{equation*}
\E(B_T| \Fc_S) \ge B_S.
\end{equation*}
\end{proof}

We note that some of the above arguments, particularly the use of
Proposition~\ref{prop:impliesmg}, allow us to deduce that if there
exists $a \in \scA$, $|a| < \infty$ for which $T \le H_a$ then
\eqref{eqn:Bdown} and \eqref{eqn:Bup} hold and $T \le H_{a'}$ for all
$a' \in \scA$.

\begin{lemma} \label{lem:old9}
Suppose that for all stopping times $S$ with $S \le T$ and
$\E|B_S|<\infty$ we have
\begin{equation} \label{eqn:Scondineq3bA}
\E(B_T | \Fc_{S}) \le B_{S} \hspace{0.4cm} \as.
\end{equation}
Then $T$ is minimal.
\end{lemma}

We refer the reader to Lemma~8 of \citet{CoxHobson:03}, the proof of which
is still valid in the more general case.

Of course we may replace the `$\le$' in \eqref{eqn:Scondineq3bA} with
`$\ge$' or `$=$' without altering the conclusion.

\begin{lemma} \label{lem:impliesminimal}
If $T \le H_{\scA} = \inf \{ t \ge 0 : B_t \in \scA\}$ is a stopping
time of the Brownian motion $(B_t)_{t \ge 0}$ where $B_0 \sim \mu_0$
and $B_T \sim \mu$, and
\begin{eqnarray}
\E(B_T | \Fc_S) & \le & B_S : \mbox{ on } \{ B_0 \ge a_- \} \label{eqn:le}\\
\E(B_T | \Fc_S) & \ge & B_S : \mbox{ on } \{ B_0 \le a_+ \}, \label{eqn:ge}
\end{eqnarray}
where $a_- = \inf \scA$ and $a_+ = \sup \scA$, then $T$ is minimal.
\end{lemma}

\begin{proof}
Choose $a \in \scA$. By assumption $T \le H_a$ and by
Lemma~\ref{lem:old9} $T$ is minimal for $\check{\mu}^a$ on $\{B_0 \le
a\}$ and for $\hat{\mu}^a$ on $\{B_0 \ge a\}$. It must then be minimal
for $\mu$.
\end{proof}

These results show the equivalence of minimality and the conditions in
\eqref{eqn:le}, \eqref{eqn:ge}. The following theorem states this
together with some extra equivalent conditions. It should be thought
of as the extension of Theorem~\ref{thm:mainB} to the setting with a
general starting measure.

\begin{theorem} \label{thm:TFAE}
Let $B$ be a Brownian motion such that $B_0 \sim \mu_0$ and $T$ a
stopping time such that $B_T \sim \mu$, where $\mu_0, \mu$ are
integrable. Let $\scA$ be the set defined in \eqref{eqn:scAdefn} and
$a_+ = \sup\{x \in [-\infty,\infty]: x \in\scA\}$, $a_- = \inf\{x \in
[-\infty,\infty]: x \in\scA\}$. Then the following are equivalent:
\begin{enumerate}
\item
$T$ is minimal;
\item
$T \le H_{\scA}$ and for all stopping times $R \le S \le T$
\begin{eqnarray*}
\E(B_S | \Fc_R) & \le & B_R \mbox{ on } \{B_0 \ge a_-\}\\
\E(B_S | \Fc_R) & \ge & B_R \mbox{ on } \{B_0 \le a_+\};
\end{eqnarray*}
\item
$T \le H_{\scA}$ and for all stopping times $S \le T$
\begin{eqnarray*}
\E(B_T | \Fc_S) & \le & B_S \mbox{ on } \{B_0 \ge a_-\}\\
\E(B_T | \Fc_S) & \ge & B_S \mbox{ on } \{B_0 \le a_+\};
\end{eqnarray*}
\item
$T \le H_{\scA}$ and for all $\gamma > 0$
\begin{eqnarray*}
\E(B_T;T > H_{-\gamma}, B_0 \ge a_-) & \le & -\gamma \Pr(T > H_{-\gamma}, B_0 \ge a_-)\\
\E(B_T;T > H_{\gamma}, B_0 \le a_+) & \ge & \gamma \Pr(T > H_{\gamma}, B_0 \le a_+);
\end{eqnarray*}
\item
$T \le H_{\scA}$ and as $\gamma \to \infty$
\begin{eqnarray*}
\gamma \Pr(T > H_{-\gamma}, B_0 \ge a_-) & \to & 0\\
\gamma \Pr(T > H_{\gamma}, B_0 \le a_+) & \to & 0.
\end{eqnarray*}
\end{enumerate}
\end{theorem}

We begin by proving the following result:

\begin{prop} \label{prop:Smodfinite}
If {\it(v)} holds and $S \le T$ then $\E|B_S| < \infty$.
\end{prop}

\begin{proof}
We show that $\E(|B_S|; B_0 \ge a_-) < \infty$. Since $B_{t \wedge
H_{-k}}$ is a supermartingale on $\{B_0 \ge -k\}$,
\begin{equation*}
\begin{split}
\E(B_{T\wedge H_{-k}};B_S < 0, & S<H_{-k}, B_0 \ge a_- \wedge (-k)) 
\\& \le \E(B_{S \wedge H_{-k}} ; B_S < 0, S < H_{-k}, B_0 \ge a_-\wedge (-k)).
\end{split}
\end{equation*}
The term on the left hand side is equal to:
\begin{equation*}
\begin{split}
\E(B_{T};B_S < 0, & T<H_{-k}, B_0 \ge a_-\wedge (-k))\\&{} - k \Pr(B_S < 0, S 
 \le H_{-k} < T, B_0 \ge a_-\wedge (-k)).
\end{split}
\end{equation*}
The first term converges (by dominated convergence) to $\E(B_T;
B_S <0, B_0 \ge a_-)$ and the second term vanishes by the
assumption. By monotone convergence
\begin{eqnarray*}
\E(B_S;B_S < 0, B_0 \ge a_-) & = & \lim_k \E(B_S; B_S < 0, S < H_{-k}, B_0 \ge a_-\wedge (-k))\\
& \ge & \lim_k \E(B_T; B_S < 0, S < H_{-k}, B_0 \ge a_-\wedge (-k))\\
& \ge & \E(B_T; B_S < 0,B_0 \ge a_-) \ge - \E(B_T^-) > -\infty.
\end{eqnarray*}

Also
\begin{eqnarray*}
\E(B_0; B_0 \ge a_- \wedge (-k)) & \ge & \E(B_{S \wedge H_{-k}}; B_0 \ge a_- \wedge(-k))\\
& = & \E(B_S;B_0 \ge a_- \wedge(-k), S<H_{-k}) \\
&&{}- k \Pr(H_{-k} \le S, B_0 \ge a_- \wedge (-k)),
\end{eqnarray*}
and
\begin{eqnarray*}
\E(B_S;B_0 \ge a_- \wedge(-k), S<H_{-k}) & = & \E(B_S^+;B_0 \ge a_-
\wedge(-k), S<H_{-k}) \\&&{}- \E(B_S^-;B_0 \ge a_- \wedge(-k),
S<H_{-k}),
\end{eqnarray*}
so
\begin{eqnarray*}
\E(B_S^+;B_0 \ge a_- \wedge(-k), S<H_{-k}) & \le & \E(B_0; B_0 \ge a_-
\wedge (-k))\\&&{} + \E(B_S^-;B_0 \ge a_- \wedge(-k), S<H_{-k})\\&&{} +
k \Pr(H_{-k} \le S, B_0 \ge a_- \wedge (-k)).
\end{eqnarray*}
By monotone and dominated convergence, in the limit we have
\begin{eqnarray*}
\E(B_S^+; B_0 \ge a_-) & \le & \E(B_0 ; B_0 \ge a_-) + \E(B_S^-; B_0 \ge a_-)\\
& < & \infty.
\end{eqnarray*}
So $\E(|B_S|; B_0 \ge a_-) < \infty$. Similarly $\E(|B_S|; B_0 \le
a_+) < \infty$, and together these imply $\E(B_S) < \infty$.
\end{proof}

\begin{proof}[Proof of Theorem~\ref{thm:TFAE}.]
Clearly {\it(ii)} $\implies$ {\it(iii)} $\implies$ {\it(iv)} $\implies$ {\it(v)} (the final
implication following from dominated convergence). We also know {\it(i)}
$\iff$ {\it(iii)}. We show {\it(v)} $\implies$ {\it(ii)}.

Suppose $A \in \Fc_R$, $A \subseteq \{ B_0 \ge a_-\}$ and set $A_k = A
\cap \{R < H_{-k}\}\cap \{B_0 \ge -k\}$. Then
\begin{equation*}
\E(B_{S \wedge H_{-k}} ; A_k) \le \E(B_{R \wedge H_{-k}};A_k).
\end{equation*}
By Proposition~\ref{prop:Smodfinite} we may apply dominated
convergence to deduce that in the limit as $k \to \infty$ the
right-hand side converges to $\E(B_R;A)$. Also
\begin{equation*}
\begin{split}
\E(B_{S \wedge H_{-k}} ; A_k) = &\E(B_S;A \cap \{B_0 \ge -k\} \cap \{S \le H_{-k}\}) \\&{}+ k \Pr(A, R < H_{-k} < S, B_0 \ge -k),
\end{split}
\end{equation*}
where the second term converges to zero by assumption and the first
converges to $\E(B_S;A)$ by dominated convergence.
\end{proof}

%****************************************************************
%****************************************************************
%****************************************************************
%****************************************************************

%% file: minlimit.tex
\section{Minimality of the Limit} \label{sec:minlimit}

We will want to show that stopping times constructed using the
techniques of Section~\ref{sec:basicbalayage} are indeed minimal when \eqref{eqn:uccond} is
satisfied. To deduce that a stopping time $T$ constructed using the
balayage techniques is minimal, we approximate $T$ by the sequence of
stopping times $T_n$ given in the construction (so $T_1$ is the exit
time from the first interval we construct, and so on). Then it is
clear that the stopping times $T_n$ satisfy the conditions of
Lemma~\ref{lem:impliesminimal}, since they are either the first exit
time from a bounded interval, or the first time to leave
$(-\infty,\alpha]$ (resp. $[\beta,\infty)$) for some $\alpha<a_-$ (resp. $\beta>a_+$). Our aim is then to deduce that
the limit is minimal.
% We shall do this by extending
%Proposition~\ref{prop:oldlimitA} to the case of a general starting
%measure.

\begin{prop} \label{prop:oldlimit}
Suppose that $T_n$ embeds $\mu_n$, $\mu_n$ converges weakly to $\mu$
and $\Pr(|T_n - T| > \eps) \to 0$ for all $\eps>0$. Then $T$
embeds $\mu$.

If also $l_n \rightarrow l_\infty < \infty$ where $l_n = \int |x|
\mu_n(dx)$ and $l_\infty = \int |x| \mu(dx)$, and $T_n$ is minimal for
$\mu_n$, then $T$ is minimal for $\mu$.
\end{prop}

\begin{remark} \label{rem:Scheffe}
Since $\mu_n \implies \mu$, on some probability space we are able to
find random variables $X_n$ and $X$ with laws $\mu_n$ and $\mu$ such
that $X_n \to X$ \as{}. By Scheff\'{e}'s Lemma therefore
\begin{equation*}
\E|X_n - X| \to 0 \mbox{ if and only if } \E|X_n| \to \E|X|,
\end{equation*}
the second statement being equivalent to $l_n \to l_{\infty}$ in the
statement of Proposition~\ref{prop:oldlimit}.
\end{remark}

Before we prove this result, we will show a useful result on the
distribution of the maximum. This will be used in the proof of Proposition~\ref{prop:oldlimit}, and also be important for the work in the next section,
when we will show that the inequality in \eqref{eqn:maxineqC} can be
attained by a class of stopping times created by balayage
techniques.

\begin{lemma} \label{lem:maxbound}
Let $T$ be a minimal embedding of $\mu$ in a Brownian motion started
with distribution $\mu_0$. Then for all $x \in \R$
\begin{equation} \label{eqn:maxineqC}
\Pr(\ol{B}_T \ge x) \le \inf_{\lambda<x} \half \left[ 1 +
\frac{u_{\mu_0}(x) - c(\lambda)}{x-\lambda} \right].
\end{equation}
\end{lemma}

\begin{proof}
Define the stopping time $\bar{H}_x = \inf\{t \ge 0 : B_t \ge x\}$, the first time that $B_t$ goes above $x$. Then we note the following inequality, which (by considering on a case by
case basis) holds for all paths and all pairs $\lambda<x$:
\begin{equation} \label{eqn:omegaineq}
\indic{\ol{B}_T \ge x} \le \frac{1}{x-\lambda}\left[ B_{T \wedge \bar{H}_x} +
\frac{|B_T - \lambda| - (B_T+\lambda)}{2} - \frac{|B_0-x| +
(B_0-x)}{2} \right].
\end{equation}
In particular, on $\{\ol{B}_T < x\}$, when therefore $\{B_0 < x\}$:
\begin{equation} \label{eqn:zerobd}
0 \le \frac{1}{x-\lambda} \left[ B_T +
\begin{Bmatrix}
-\lambda & : B_T > \lambda \\ -B_T & : B_T \le \lambda
\end{Bmatrix}
\right].
\end{equation}
While on $\{\ol{B}_T \ge x\}$,
\begin{eqnarray}
1 & \le & \frac{1}{x-\lambda}\left[ B_{T \wedge \bar{H}_x} +
  \begin{Bmatrix}
    -\lambda & : B_T > \lambda \\ -B_T & : B_T \le \lambda
  \end{Bmatrix}
-
  \begin{Bmatrix}
  B_0 -x & : B_0 > x  \\ 0 & : B_0 \le x
  \end{Bmatrix}
\right]\nonumber\\
& \le & \frac{1}{x-\lambda}\left[ x + 
  \begin{Bmatrix}
    -\lambda & : B_T > \lambda \\ -B_T & : B_T \le \lambda
  \end{Bmatrix}
\right].\label{eqn:onebd}
\end{eqnarray}

So we may take expectations in \eqref{eqn:omegaineq} to get
\begin{equation} \label{eqn:maxforlambda}
\Pr(\ol{B}_T \ge x) \le \frac{1}{2} \left[ 1 + \frac{2\E(B_{T \wedge \bar{H}_x})
+ \left(u_{\mu_0}(x) - u_\mu(\lambda)\right) - \left(\E(B_T) +
\E(B_0)\right)}{(x-\lambda)} \right].
\end{equation}
We can deduce \eqref{eqn:maxineqC} provided we can show
\begin{equation} \label{eqn:CgeTwedgeHx}
C \ge 2\E(B_{T \wedge \bar{H}_x})-\left(\E(B_T) + \E(B_0)\right)
\end{equation}
since \eqref{eqn:maxforlambda} holds for all $\lambda<x$.

We now consider $a \in \scA$ possibly taking the values $\pm
\infty$. Since $u_{\mu}(a) - u_{\mu_0}(a) = C$ for $a \in \scA$, we
can deduce
\begin{equation*}
C = 2\E(B_T;B_T \ge a) + 2 \E(B_0;B_0 < a) - \E(B_T) - \E(B_0)
\end{equation*} 
where we note that $\{B_T < a\} = \{B_0 < a\}$. Theorem~\ref{thm:TFAE}
tells us that
\begin{eqnarray}
\E(B_{T \wedge \bar{H}_x};B_0 < a) & \le &  \E(B_T;B_0 < a)\label{eqn:EBTHxle1}\\
\E(B_{T \wedge \bar{H}_x};B_0 \ge a) & \le & \E(B_0;B_0 \ge a)\label{eqn:EBTHxle2}
\end{eqnarray}
\begin{comment}
\begin{equation}\label{eqn:EBTHxle}
\E(B_{T \wedge H_x}) \le \E(B_T;B_0 < a) + \E(B_0;B_0 \ge a)
\end{equation}
\end{comment}
and \eqref{eqn:CgeTwedgeHx} holds.
\end{proof}

We also have the following result:
\begin{prop} \label{prop:unifconv}
Suppose $\mu$ and $\{\mu_n\}_{n \ge 1}$ are all integrable
distributions such that $\mu_n \implies \mu$ and $l_n =
\int |y| \, \mu_n(dy) \to \int |y| \, \mu(dy) = l_\infty$. Then
$u_{\mu_n}$ converges uniformly to $u_{\mu}$.
\end{prop}
\begin{proof}
Fix $\eps > 0$. By \eqref{eqn:umuis}, using the fact that $\mu - \hat
\mu = \check \mu$ we may write
\begin{equation*}
u_\mu (x) = \int_{-\infty}^{\infty} (x-y) \, \mu(dy) + 2
\int_{-\infty}^x (y-x) \, \mu(dy) = x - \int_{-\infty}^\infty y \,
\mu(dy) + 2 \int_{-\infty}^x (y-x) \, \mu(dy),
\end{equation*}
and similarly for $u_{\mu_n}$, hence
\begin{equation} \label{eqn:umuminumun}
u_{\mu_n}(x) - u_{\mu}(x) = (m_\infty - m_n) + 2 \int_{-\infty}^x(y-x) \, (\mu_n - \mu)(dy),
\end{equation}
where we write $m_n, m_\infty$ for the means of $\mu_n$ and $\mu$
respectively; $m_n \to m_\infty$ as a consequence of Remark~\ref{rem:Scheffe}.  Since $\mu$ is integrable, as $x \downarrow -\infty$,
\begin{equation*}
\int_{-\infty}^x (x-y) \, \mu(dy) \downarrow 0.
\end{equation*}
By \eqref{eqn:umuminumun} and Lemma~\ref{lem:limexists} (which implies
$u_{\mu_n}$ converges to $u_\mu$ pointwise, the $C$ in \eqref{eqn:munlimit} being $0$ since $l_n \to l_\infty$), for all $x \in \R$
\begin{equation*}
\int_{-\infty}^x (x-y) \mu_n(dy) \to \int_{-\infty}^x (x-y) \mu(dy)
\end{equation*}
as $n \to \infty$. Finally we note that both sides of the above are
increasing in $x$.

Consider
\begin{equation*}
|u_{\mu_n}(x) - u_{\mu}(x)| \le |m_\infty - m_n| + 2 \int_{-\infty}^x (x-y) \, \mu_n(dy) + 2 \int_{-\infty}^x (x-y)\, \mu(dy).
\end{equation*}
We may choose $x_0$ sufficiently small that $\int_{-\infty}^{x_0} (x_0
- y) \, \mu(dy) < \eps$, and therefore such that
\begin{equation*}
\int_{-\infty}^{x} (x - y) \, \mu(dy) \le \int_{-\infty}^{x_0} (x_0 - y) \, \mu(dy) < \eps
\end{equation*}
for all $x \le x_0$. By the above and Remark~\ref{rem:Scheffe} we may
now choose $n_0(\eps)$ such that for all $n \ge n_0(\eps)$
\begin{equation*}
|m_\infty - m_n| < \eps \mbox{ and } \left| \int_{-\infty}^{x_0}(x_0 - y) \, \mu_n(dy) - \int_{-\infty}^{x_0}(x_0 - y) \, \mu(dy) \right| < \eps.
\end{equation*}
Then for all $x \le x_0$ and for all $n \ge n_0(\eps)$,
\begin{equation*}
|u_{\mu_n}(x) - u_{\mu}(x)| \le \eps + 2 \times 2 \eps + 2 \eps = 7 \eps.
\end{equation*}

Similarly we can find $x_1, n_1(\eps)$ such that $|u_{\mu_n}(x) -
u_{\mu}(x)| \le 7 \eps$ for all $x \ge x_1$ and all $n \ge
n_1(\eps)$. Finally $u_{\mu_n}, u_\mu$ are both Lipschitz and
pointwise $u_{\mu_n}(x) \to u_{\mu}(x)$ and we must have uniform
convergence on any bounded interval, and in particular on $[x_0,x_1]$.
\end{proof}

\begin{proof}[Proof of Proposition~\ref{prop:oldlimit}.]
Suppose first that there exists $a \in \scA \cap \R$. We show that $T
\le H_a$ for all such $a$. As usual, we write $\mu_0$ for the starting
measure, and $c(x)=u_{\mu}(x) - C$. We define $C_n$ to be the smallest
value such that $u_{\mu_0}(x) \ge u_{\mu_n}(x) -C_n$ and the functions
$c_n(x) = u_{\mu_n}(x) -C_n$. Note that $l_n = u_{\mu_n}(0)$, so
$\lim_{n \to \infty} u_{\mu_n}(0)$ exists. Then
(by Lemma~\ref{lem:limexists}{\it (i)} or equivalently \citep{Chacon:77}[Lemma 2.5]) weak
convergence implies
\begin{equation*}
\lim_{n \to \infty} u_{\mu_n}(x) = u_\mu(x) - K
\end{equation*}
for all $x \in \R$ and (here) $K = 0$ since $u_{\mu_n}(0) \to
u_{\mu}(0)$.
%
%\begin{comment}
% We note also that weak convergence together with $l_n \to l_\infty$
%implies that $\int x \, \mu_n(dx) \to \int x \, \mu(dx)$. Recall
%\begin{equation*}
%C = \inf_{x \in \R}\{u_\mu(x) - u_{\mu_0}(x)\}.
%\end{equation*}
%Suppose this is attained at $x_0 \in \R$ (the case where the limit is
%attained as $x \to \pm \infty$ follows from
%\eqref{eqn:pmulimit}). Then since as $n \to \infty$, $u_{\mu_n}(x_0)
%\to u_{\mu}(x_0)$ and
%\begin{equation*}
%C_n \le u_{\mu_n}(x_0) - u_{\mu_0}(x_0),
%\end{equation*}
%give $\eps > 0$ we may find $n_0(\eps)$ such that $n\ge n_0$ implies
%$C_n \le C+\eps$. Then for such $n$,
%\begin{eqnarray*}
%\Pr(M_{T_n} \ge x) & \le & \frac{1}{2} \left[ 1 + \frac{u_{\mu_0}(x) -
%u_{\mu_n}(\lambda) + C}{x-\lambda} + \frac{C_n-C}{x- \lambda}\right]\\
%& \le & \frac{1}{2} \left[ 1 + \frac{u_{\mu_0}(x) - u_{\mu_n}(\lambda)
%+ C}{x-\lambda} + \frac{\eps}{x- \lambda}\right].
%\end{eqnarray*}
%We take the limit as $n \to \infty$, noting that $\Pr(M_{T_n} \ge x)
%\to \Pr(M_T \ge x)$ and
%\begin{eqnarray*}
%\Pr(M_{T} \ge x) & \le & \frac{1}{2} \left[ 1 + \frac{u_{\mu_0}(x) -
%c(\lambda)}{x-\lambda} + \frac{\eps}{x- \lambda}\right]\\ & \le &
%\frac{1}{2} \left[ 1 + \frac{u_{\mu_0}(x) -
%c(\lambda)}{x-\lambda}\right],
%\end{eqnarray*}
%since $\eps$ is arbitrary. 
%\end{comment}

By Lemma~\ref{lem:maxbound} for $x \in \R$ and $\lambda < x$
\begin{equation*}
\Pr(\ol{B}_{T_n} \ge x)  \le  \frac{1}{2} \left[ 1 + \frac{u_{\mu_0}(x) -
u_{\mu_n}(\lambda) + C}{x-\lambda} + \frac{C_n-C}{x- \lambda}\right],
\end{equation*}
and we take the limit as $n \to \infty$, using
Proposition~\ref{prop:unifconv} (so that $C_n \to C$) and noting that
$\Pr(\ol{B}_{T_n} \ge x) \to \Pr(\ol{B}_T \ge x)$, to get
\begin{equation*}
\Pr(\ol{B}_{T} \ge x) \le \frac{1}{2} \left[ 1 + \frac{u_{\mu_0}(x) -
c(\lambda)}{x-\lambda} \right].
\end{equation*}

Suppose now $x=a$. Since the above holds for all $\lambda<a$, we may
take the limit of the right hand side as $\lambda \uparrow a$, in
which case $u_{\mu_0}(a) = c(a)$, and by Remark~\ref{rem:potprop}(v)
\begin{eqnarray*}
\Pr(\ol{B}_T \ge a) & \le & \half \left[ 1+c_-'(a)\right]\\
& \le & \half \left[ 1 + (1-2\mu((-\infty,a)))\right]\\
& \le & \mu([a,\infty)).
\end{eqnarray*}
By considering $-B_t$ we may deduce that $\Pr(\ul{B}_T \le a) \le
\mu((-\infty,a])$. Hence $\Pr(T \le H_a) = 1$, and we deduce that $T$
is minimal.

It only remains to show (by Lemma~\ref{lem:old9}) that if $\infty \in
\scA$ then
\begin{equation*}
\E(B_T | \Fc_S) \ge B_S
\end{equation*}
for all stopping times $S \le T$. The case where $-\infty \in \scA$
follows from $B_t \mapsto -B_t$. In particular, for $S \le T$ and $A
\in \Fc_S$ we need to show
\begin{equation}\label{eqn:iffcondB}
\E(B_T ; A) \ge \E(B_S ; A).
\end{equation}

In fact we need only show the above for sets $A \subseteq \{S < T\}$
since it clearly holds on $\{S = T\}$. So we can define $A_n = A \cap
\{S < T_n\}$ and therefore $\Pr(A \setminus A_n) \to 0$ as $n \to
\infty$. Also $A_n \in \Fc_{S \wedge T_n}$. By Theorem~\ref{thm:TFAE}
and the fact that the $T_n$ are minimal
\begin{eqnarray*}
\E(B_{S \wedge T_n};A_n) & \le & \E(B_{T_n}; A_n \cap \{B_0 \le a_+^n\})
                + \E(B_{S \wedge T_n}; B_0 > a_+^n)
                \\ && {}- \E(B_{S \wedge T_n}; A_n^C \cap \{B_0 > a_+^n\})\\
& \le & \E(B_{T_n}; A_n \cap \{B_0 \le a_+^n\})
                \E(B_0; B_0 > a_+^n) \\ &&{}- \E(B_{T_n}; A_n^C \cap \{B_0 > a_+^n\})\\
& \le & \E(B_{T_n}; A_n) - \E(B_{T_n}; \{B_0>a_+^n\}) + \E(B_0; B_0 > a_+^n)
\end{eqnarray*}
where $a_+^n$ is the supremum of the set $\scA_n$ (that is the
corresponding set to $\scA$ for the measures $\mu_0,\mu_n$). This is
not necessarily infinite.

So it is sufficient for us to show that
\begin{eqnarray}
\lim_n \E(B_{T_n};A_n) & = & \E(B_{T};A)\label{eqn:lim1a};\\
\lim_n \E(B_S ; A_n) & = & \E(B_S ; A) \label{eqn:lim2a},
\end{eqnarray}
and
\begin{equation}
\lim_n |\E(B_0; B_0 > a_+^n) - \E(B_{T_n}; B_0>a_+^n)| = 0 \label{eqn:lim3a}.
\end{equation}

For \eqref{eqn:lim1a}
%we may use a proof identical to that used in Proposition~\ref{prop:oldlimitA} to show 
%\eqref{eqn:lim1}.
we consider $|\E(B_{T};A) -
\E(B_{T_n};A_n)|$. Then
\begin{equation*}
  |\E(B_{T};A) - \E(B_{T_n};A_n)| \le \E(|B_{T}|; A
  \setminus A_n) + \E(|B_{T} - B_{T_n}|;A_n)
\end{equation*}
and the first term tends to zero by dominated convergence (this
follows from the assumption that $T_n$ converges to $T$ in
probability). For the second term we show $\E(|B_{T} - B_{T_n}|) \to
0$. Fix $\eps >0$. We have
\begin{equation*}
  |B_{T} - B_{T_n}| \le |B_{T_n}| - |B_{T}| + 2
  |B_{T}| \indic{|T_n - T| \ge \eps} + 2|B_{T} -
  B_{T_n}| \indic{|T_n - T| \le \eps}.
\end{equation*}
We take expectations and let $n \to \infty$. By the definition of
$\mu_n$ the first two terms cancel each other out, while the third
tends to zero by dominated convergence. For the last term, by the
(strong) Markov property
\begin{equation*}
  \E(|B_{T} - B_{T_n}|;|T_n - T| \le \eps) \le
  \E(|B_\eps|) = \sqrt{\frac{\eps}{2 \pi}}.
\end{equation*}
Consequently, in the limit, $\E(|B_{T} - B_{T_n}|;|T_n - T| \le \eps)
\to 0$ and \eqref{eqn:lim1a} holds.
We want to apply Lemma~\ref{lem:old9} so we can assume that $\E|B_S| < \infty$, and
\eqref{eqn:lim2a} follows by dominated convergence.

Finally we consider \eqref{eqn:lim3a}. Let $\theta_n =
\mu_0((-\infty,a_+^n])$.  Since $a_+^n \in \scA_n$ we have
\begin{eqnarray*}
\E(B_{0}; B_0 > a_+^n) - \E(B_{T_n}; B_0 > a_+^n) & = & \int y \, \hat\mu_0^{a_+^n,\theta_n}(dy) - \int y \, \hat\mu_n^{a_+^n,\theta_n}(dy)\\
& = & \int (y - a_+^n) \, \hat\mu_0^{a_+^n,\theta_n}(dy) - \int (y - a_+^n) \, \hat\mu^{a_+^n,\theta_n}(dy)\\
& = & \half \left[ \int y \, (\mu_0 - \mu_n)(dy) + u_{\mu_n}(a_+^n) - u_{\mu_0}(a_+^n)\right]\\
& = & \half \left[ \int y \, (\mu_0-\mu)(dy) - C_n\right],
\end{eqnarray*}
where we have used the fact that (for a general measure $\nu$)
\begin{equation*}
\int (y-x) \,\hat \nu^x(dy) = \half \left[ \int y \, \nu(dy) - u_\nu(x) - x\right].
\end{equation*}
As $n \to \infty$, since $\infty \in \scA$,
\begin{equation*}
\int y \, (\mu_0 - \mu_n)(dy) \to \int y \, (\mu_0 - \mu)(dy) = C.
\end{equation*}
So we need only show that $C_n \to C$, which follows from the uniform
convergence of $u_{\mu_n}$ to $u_\mu$ (Proposition~\ref{prop:unifconv}).
\end{proof}

%****************************************************************
%****************************************************************
%****************************************************************
%****************************************************************

%% file: AYtype.tex
\section{Tangents and Azema-Yor Type Embeddings} \label{sec:AYtype}

One of the motivations for this paper is to discuss generalisations
of the Azema-Yor family of embeddings (see
\citet{AzemaYor:79,Jacka:88}) to the
integrable starting/target measures we have discussed already.

The aim is therefore to find the embedding which maximises the law of
the maximum, $\sup_{0 \le t \le T} B_t$ (or in the more general case
$\sup_{0 \le t \le T} |B_t|$). If we look for the maximum within the
class of all embeddings there is no natural maximum embedding. For
this reason we consider the class of minimal
embeddings. Lemma~\ref{lem:maxbound} establishes that there is some
natural limit when we consider this restriction. In fact the extended
Azema-Yor embedding will attain the limit in \eqref{eqn:maxineqC}.

The idea is to use the machinery from the previous sections to show
the embeddings exist as limits of the Chacon-Walsh type embeddings of
Section~\ref{sec:basicbalayage}. It is then possible to show that the
embeddings are minimal and that they attain equality in
\eqref{eqn:maxineqC}.

\begin{theorem}
If $T$ is a stopping time as described in Lemma~\ref{lem:fnembed2},
where $C$ as described in the lemma is
\begin{equation}\label{eqn:Cis}
C = \inf_{x} \{u_\mu(x) - u_{\mu_0}(x)\},
\end{equation}
then $T$ is minimal.
\end{theorem}

\begin{proof}
Lemma~\ref{lem:fnembed2} suggests a sequence $T_n$ of stopping times
for which $T$ is the limit. We note that we can modify the definition
of $T_n$ so that $T_n'$ is specified by the functions
$f_1,f_2,\ldots,f_n,f^{-1},f^{+1}$ without altering their limit (as a
consequence of \eqref{eqn:gefn}), where $f^{-1}$ is the tangent to $g$
with gradient $-1$ and $f^{+1}$ is the tangent to $g$ with gradient
$1$. It is easy to see that this ensures that $\E(B_{T_n'}) =
\E(B_T)$ (by \eqref{eqn:pmulimit}), and also that $u_{\mu_n}(0) \to u_{\mu}(0)$ and $n \to
\infty$. Consequently the stopping times $T_n'$ and $T$ satisfy the
conditions of Proposition~\ref{prop:oldlimit}, where it is clear that
the $T_n'$ are all minimal, since each step clearly satisfies the conditions of Theorem~\ref{thm:TFAE} as a consequence of \eqref{eqn:Cis}. So $T$ is minimal.
\end{proof}

Define the function
\begin{equation} \label{eqn:Phidefn}
\Phi(x) = \argmin_{\lambda < x} \left\{ \frac{u_{\mu_0}(x) - c(\lambda)}{x-\lambda} \right\}.
\end{equation}
In the cases described by \citet{AzemaYor:79}, this is the barycentre
function. It can also be seen to agree with the function appearing in
the generalisation of the Azema-Yor stopping time to non-centred means
which appears in \citet{CoxHobson:03}. A similar function is used in
\citet{Hobson:98a} who examines the case where starting and target
means are centred and satisfy \eqref{eqn:grcondition}. $\Phi(\cdot)$
can be thought of graphically as the point (below $x$) at which there
exists a tangent to $c(\cdot)$ meeting the function $u_{\mu_0}(\cdot)$
at $x$.

\begin{lemma}
The Azema-Yor stopping time
\begin{equation} \label{eqn:AYtime}
T = \inf \{ t \ge 0 : B_t \le \Phi (\ol{B}_t) \}
\end{equation}
is minimal and attains equality in \eqref{eqn:maxineqC}.
\end{lemma}

\begin{figure}[t]
\begin{center}
\includegraphics[width=\textwidth,height=3in]{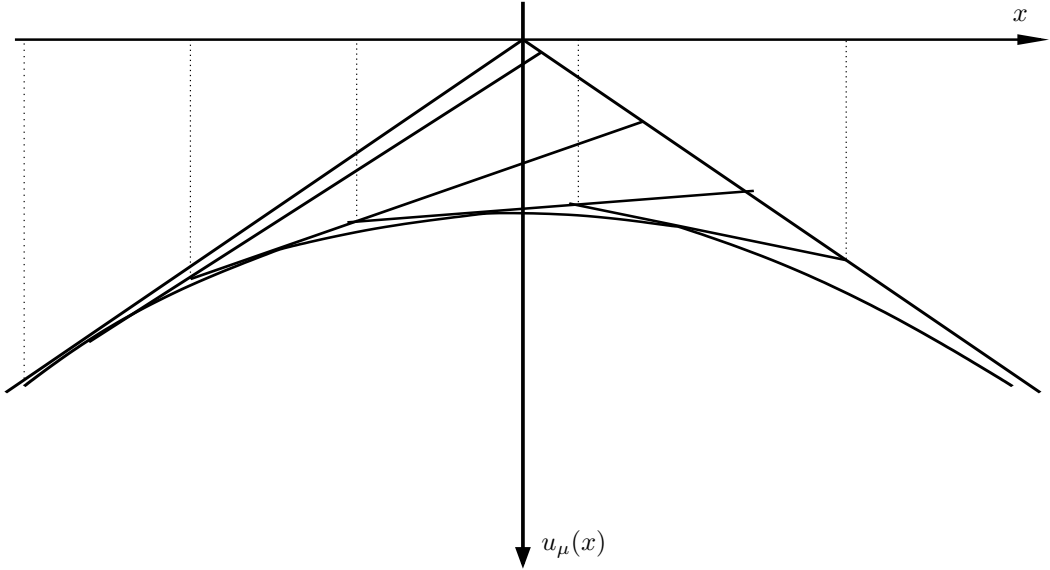}
\caption[Approximating the Azema-Yor stopping time]{\label{fig:cwfig4}Approximating the Azema-Yor stopping time: we take tangents to the potential from left to right. In the limit the tangents become closer. The dotted lines highlight the points at which the approximated stopping time will stop the process.}
\end{center}
\end{figure}

We prove this lemma using an extension of an idea first suggested in
\citet{Meilijson:83}. We approximate $T$ by taking tangents to $c$,
starting with gradient $-1$, and increasing to $+1$. As the number of
tangents we take increases, the stopping time converges to $T$.  The
general approximation sequence can be seen in Figure~\ref{fig:cwfig4}.

\psfrag{x}{$x$}
\psfrag{umu}{$u_{\mu}(x)$}

\begin{proof}
We apply Lemma~\ref{lem:fnembed2} for each $n$ to the functions
$f_1^n,f_2^n,\ldots,f_{m(n)}^n$, which are chosen as tangents to
$c(\cdot)$ with increasing gradients, so that $f_1^n$ has gradient
$-1$, $f_m^n$ has gradient $1$, and so that the difference in the
gradients of consequential tangents is less than $\frac{1}{n}$. We
also choose the tangents in such a way that the points at which successive
tangents intersect each other (which are $B_{T_n}$ stops) are at most
$\frac{1}{n}$ apart when they lie within $[-n,n]$ (at least as far as
this is possible --- if both $\mu_0$ and $\mu$ have an interval
containing no mass, it might not be possible to manage this, but this
case will not be important). This defines a (minimal) stopping time
$T_n$ such that (by \eqref{eqn:pmulimit}) $\E(B_{T_n}) = \int x \,
\mu(dx)$. Also, by considering $\mu_n = \Lc(B_{T_n})$,
$|\mu_n((-\infty,x)) - \mu((-\infty,x))| \le \frac{1}{n}$ for all $x
\in \R$. So $\mu_n \implies \mu$. The choice of $T_n$ also ensures
that $\Pr(|T-T_n| > \eps) \to 0$ for all $\eps > 0$. Consequently $T$
is minimal.

To deduce that $T$ attains equality in \eqref{eqn:maxineqC} we note
that $\Phi(x)$ is the optimal choice for $\lambda$ in
\eqref{eqn:maxineqC}, and by the definition of $\Phi(x)$,
\begin{eqnarray*}
\{\ol{B}_T < x \} & \subseteq & \{B_T \le \Phi(x)\}\\
\{\ol{B}_T \ge x \} & \subseteq & \{B_T \ge \Phi(x)\}.
\end{eqnarray*}
This means we attain equality in \eqref{eqn:zerobd} and
\eqref{eqn:onebd}, and so only need show that we have equality in
\eqref{eqn:EBTHxle1} and \eqref{eqn:EBTHxle2} for equality in
\eqref{eqn:maxineqC} to hold. But for $x$ given, we may calculate the
potential of $\mu' = \Lc(B_{T \wedge \bar{H}_x})$ --- where $\bar{H}_x =
\inf\{t \ge 0: B \ge x\}$ --- as:
\begin{equation*}
u_{\mu'}(y) = \begin{cases}
u_{\mu}(y) &: y \le \Phi(x);\\
u_{\mu}(\Phi(x)) + \frac{y-\Phi(x)}{x-\Phi(x)}(u_{\mu_0}(x) - u_{\mu}(\Phi(x))) & : \Phi(x) \le y \le x;\\
u_{\mu_0}(y) &: y \ge x.
\end{cases}
\end{equation*}
It then follows from Theorem~\ref{thm:TFAE} and \eqref{eqn:pmulimit}
that equality holds.
\end{proof}

\begin{remark}
The embedding due to \citet{Jacka:88} can be viewed easily in this framework. Essentially the embedding can be described as follows. We wish to find an embedding which maximises $\Pr(\sup_{t \le T} |B_t| \ge x)$ simultaneously for all $x$, subject to $T$ being minimal. The optimal construction can be described thus in the Chacon-Walsh picture: choose the tangent to $u_{\mu}$  with gradient 0. Let $(a,b)$ be the points at which the tangent intersects $c(x)$, where $C$ is chosen to be the value given by \eqref{eqn:uccond}; the first step of the construction is to run the process until it leaves the interval $I=(a,b)$. The construction will now have at least two separate halves; on the positive half we run the Azema-Yor construction of \eqref{eqn:AYtime}, and on the negative half we run the reverse of the Azema-Yor construction. Such a construction will therefore be minimal, and can be shown to be optimal --- for details we refer the reader to \citet{CoxHobson:03}.

This technique can be trivially extended to maximising $\Pr(\sup_{t \le T} f(B_t) \ge x)$, where $f$ is a function which is increasing above some point $x_0$ and decreasing below $x_0$.
\end{remark}

\psfrag{e}{$\eps$}
\begin{figure}[t]
\begin{center}
\includegraphics[width=\textwidth]{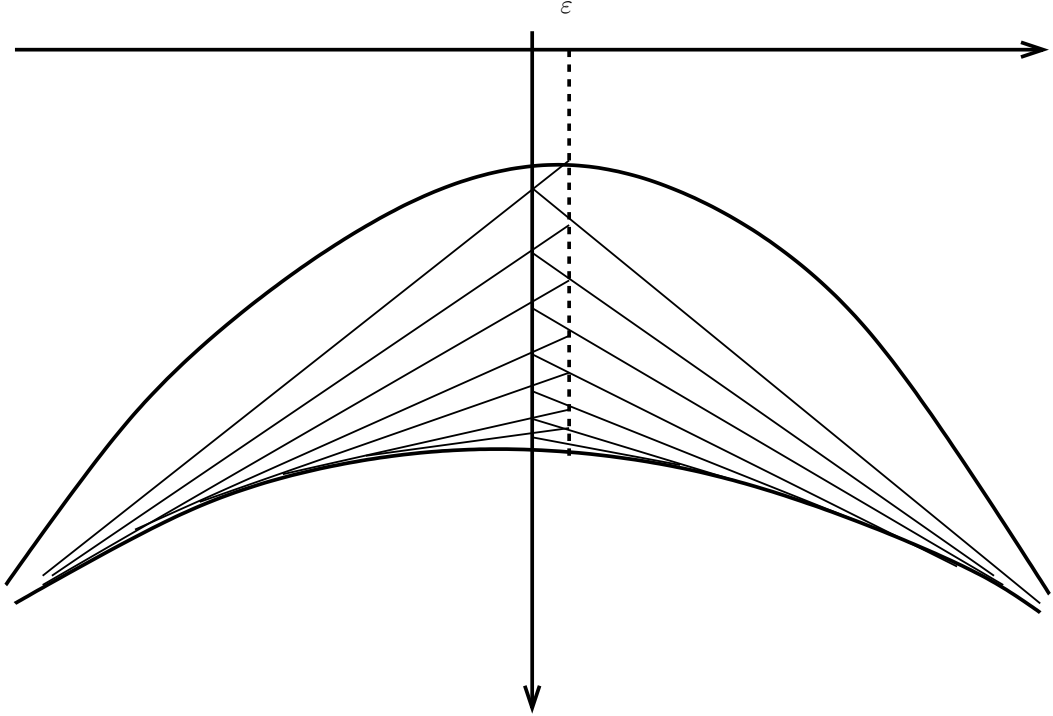}
\caption[The Chacon-Walsh type picture for an
approximation to the Vallois stopping time, with a general starting
distribution]{\label{fig:vallois}The Chacon-Walsh type picture for an
approximation to the Vallois stopping time, with a general starting
distribution. We note that after the first two steps, there could
still be mass at the extremes. This mass will have to be embedded
using some suitable procedure --- for example a Vallois construction
using the local time at a different level.}
\end{center}
\end{figure}

\begin{remark}
In a slightly different vein, the construction of \citet{Vallois:83} with decreasing functions can also be seen in this framework. Choose $\eps>0$ and construct alternate tangents; tangents of positive gradient intersecting the current potential at $\eps$, and tangents of negative gradient intersecting the current potential at 0 (see Figure~\ref{fig:vallois}). This can be repeated a number of times, and results in an approximation of $c(x)$; suitable further choices of tangents can be used to construct a full embedding. If $c$ touches $u_{\mu_0}$ at finite points, any choice of construction outside the interval containing $0$ may be used. By construction this is a minimal embedding.

The Vallois construction results on taking the limit as $\eps \downarrow 0$; the appearance of the local time in the construction results since the limiting embeddings are determined by the number of downcrossings of $[0,\eps]$, which has as a limit the local time (see \citet{RevuzYor:99}[VI.1.10]). An application of Proposition~\ref{prop:oldlimit} allows us to deduce that this limit is minimal. For further details we refer the reader to \citet{Cox:04}.
\end{remark}